\documentclass[a4,12pt]{article}
\usepackage{amssymb}
\input diagrams
\input prepictex
\input pictex
\input postpictex

\begin{document}

\baselineskip15.5pt

\newtheorem{de}{Definition $\!\!$}[section]
\newtheorem{pr}[de]{Proposition $\!\!$}
\newtheorem{lem}[de]{Lemma $\!\!$}
\newtheorem{corollary}[de]{Corollary $\!\!$}
\newtheorem{them}[de]{Theorem $\!\!$}
\newtheorem{example}[de]{\it Example $\!\!$}
\newtheorem{remark}[de]{\it Remark $\!\!$}

\newcommand{\nc}[2]{\newcommand{#1}{#2}}
\newcommand{\rnc}[2]{\renewcommand{#1}{#2}}

\rnc{\theequation}{\thesection.\arabic{equation}}

\def\Box{\square}
\def\Diamond{\lozenge}
\def\text{\mbox}
\def\eps{{\epsilon}}
\nc{\beq}{\begin{equation}}
\nc{\eeq}{\end{equation}}
\rnc{\[}{\beq}
\rnc{\]}{\eeq}
\nc{\qpb}{quantum principal bundle}
\nc{\bpr}{\begin{pr}}
\nc{\bth}{\begin{them}}
\nc{\ble}{\begin{lem}}
\nc{\bco}{\begin{corollary}}
\nc{\bre}{\begin{remark}}
\nc{\bex}{\begin{example}}
\nc{\bde}{\begin{de}}
\nc{\ede}{\end{de}}
\nc{\epr}{\end{pr}}
\nc{\ethe}{\end{them}}
\nc{\ele}{\end{lem}}
\nc{\eco}{\end{corollary}}
\nc{\ere}{\hfill\mbox{$\Diamond$}\end{remark}}
\nc{\eex}{\hfill\mbox{$\Diamond$}\end{example}}
\nc{\epf}{\hfill\mbox{$\Box$}\vspace*{3mm}}
\nc{\ot}{\otimes}
\nc{\bsb}{\begin{Sb}}
\nc{\esb}{\end{Sb}}
\nc{\ct}{\mbox{${\cal T}$}}
\nc{\ctb}{\mbox{${\cal T}\sb B$}}
\nc{\ba}{\begin{array}}
\nc{\ea}{\end{array}}
\nc{\bea}{\begin{eqnarray}}
\nc{\beas}{\begin{eqnarray*}}
\nc{\eeas}{\end{eqnarray*}}
\nc{\eea}{\end{eqnarray}}
\nc{\be}{\begin{enumerate}}
\nc{\ee}{\end{enumerate}}
\nc{\bcd}{\beq\begin{CD}}
\nc{\ecd}{\end{CD}\eeq}
\nc{\bi}{\begin{itemize}}
\nc{\ei}{\end{itemize}}
\nc{\kr}{\mbox{Ker}}
\nc{\te}{\!\ot\!}
\nc{\pf}{\mbox{$P\!\sb F$}}
\nc{\pn}{\mbox{$P\!\sb\nu$}}
\nc{\bmlp}{\mbox{\boldmath$\left(\right.$}}
\nc{\bmrp}{\mbox{\boldmath$\left.\right)$}}
\rnc{\phi}{\varphi}
\nc{\LAblp}{\mbox{\LARGE\boldmath$($}}
\nc{\LAbrp}{\mbox{\LARGE\boldmath$)$}}
\nc{\Lblp}{\mbox{\Large\boldmath$($}}
\nc{\Lbrp}{\mbox{\Large\boldmath$)$}}
\nc{\lblp}{\mbox{\large\boldmath$($}}
\nc{\lbrp}{\mbox{\large\boldmath$)$}}
\nc{\blp}{\mbox{\boldmath$($}}
\nc{\brp}{\mbox{\boldmath$)$}}
\nc{\LAlp}{\mbox{\LARGE $($}}
\nc{\LArp}{\mbox{\LARGE $)$}}
\nc{\Llp}{\mbox{\Large $($}}
\nc{\Lrp}{\mbox{\Large $)$}}
\nc{\llp}{\mbox{\large $($}}
\nc{\lrp}{\mbox{\large $)$}}
\nc{\lbc}{\mbox{\Large\boldmath$,$}}
\nc{\lc}{\mbox{\Large$,$}}
\nc{\Lall}{\mbox{\Large$\forall\;$}}
\nc{\bc}{\mbox{\boldmath$,$}}
\rnc{\epsilon}{\varepsilon}
\nc{\ra}{\rightarrow}
\nc{\ci}{\circ}
\nc{\cc}{\!\ci\!}
\nc{\T}{\mbox{\sf T}}
\nc{\can}{\mbox{\em\sf T}\!\sb R}
\nc{\cnl}{$\mbox{\sf T}\!\sb R$}
\nc{\lra}{\longrightarrow}
\nc{\M}{\mbox{\rm Map}}
\rnc{\to}{\mapsto}
\rnc{\breve}{\widetilde}
\nc{\imp}{\Rightarrow}
\rnc{\iff}{\Leftrightarrow}
\nc{\ob}{\mbox{$\Omega\sp{1}\! (\! B)$}}
\nc{\op}{\mbox{$\Omega\sp{1}\! (\! P)$}}
\nc{\oa}{\mbox{$\Omega\sp{1}\! (\! A)$}}
\nc{\inc}{\mbox{$\,\subseteq\;$}}
\rnc{\subset}{\inc}
\nc{\spp}{\mbox{${\cal S}{\cal P}(P)$}}
\nc{\dr}{\mbox{$\Delta_{R}$}}
\nc{\dsr}{\mbox{$\Delta_{\Omega^1P}$}}
\nc{\ad}{\mbox{$\mathop{\mbox{\rm Ad}}_R$}}
\nc{\m}{\mbox{m}}
\nc{\0}{\sb{(0)}}
\nc{\1}{\sb{(1)}}
\nc{\2}{\sb{(2)}}
\nc{\3}{\sb{(3)}}
\nc{\4}{\sb{(4)}}
\nc{\5}{\sb{(5)}}
\nc{\6}{\sb{(6)}}
\nc{\7}{\sb{(7)}}
\nc{\hsp}{\hspace*}
\nc{\nin}{\mbox{$n\in\{ 0\}\!\cup\!{\Bbb N}$}}
\nc{\al}{\mbox{$\alpha$}}
\nc{\bet}{\mbox{$\beta$}}
\nc{\ha}{\mbox{$\alpha$}}
\nc{\hb}{\mbox{$\beta$}}
\nc{\hg}{\mbox{$\gamma$}}
\nc{\hd}{\mbox{$\delta$}}
\nc{\he}{\mbox{$\varepsilon$}}
\nc{\hz}{\mbox{$\zeta$}}
\nc{\hs}{\mbox{$\sigma$}}
\nc{\hk}{\mbox{$\kappa$}}
\nc{\hm}{\mbox{$\mu$}}
\nc{\hn}{\mbox{$\nu$}}
\nc{\la}{\mbox{$\lambda$}}
\nc{\hl}{\mbox{$\lambda$}}
\nc{\hG}{\mbox{$\Gamma$}}
\nc{\hD}{\mbox{$\Delta$}}
\nc{\Th}{\mbox{$\Theta$}}
\nc{\ho}{\mbox{$\omega$}}
\nc{\hO}{\mbox{$\Omega$}}
\nc{\hp}{\mbox{$\pi$}}
\nc{\hP}{\mbox{$\Pi$}}
\nc{\bpf}{{\it Proof.~~}}
\nc{\as}{\mbox{$A(S^3\sb s)$}}
\nc{\bs}{\mbox{$A(S^2\sb s)$}}
\nc{\slq}{\mbox{$A(SL\sb q(2))$}}
\nc{\fr}{\mbox{$Fr\llp A(SL(2,\IC))\lrp$}}
\nc{\slc}{\mbox{$A(SL(2,\IC))$}}
\nc{\af}{\mbox{$A(F)$}}
\nc{\suq}{\mbox{$A(SU_q(2))$}}
\nc{\asq}{\mbox{$A(S_q^2)$}}
\nc{\tasq}{\mbox{$\widetilde{A}(S_q^2)$}}
\nc{\tc}{\widetilde{can}}

\def\esl{{\mbox{$E\sb{\frak s\frak l (2,{\Bbb C})}$}}}
\def\esu{{\mbox{$E\sb{\frak s\frak u(2)}$}}}
\def\zf{{\mbox{${\Bbb Z}\sb 4$}}}
\def\zt{{\mbox{$2{\Bbb Z}\sb 2$}}}
\def\ox{{\mbox{$\Omega\sp 1\sb{\frak M}X$}}}
\def\oxh{{\mbox{$\Omega\sp 1\sb{\frak M-hor}X$}}}
\def\oxs{{\mbox{$\Omega\sp 1\sb{\frak M-shor}X$}}}
\def\Fr{\mbox{Fr}}
\def\gal{-Galois extension}
\def\hge{Hopf-Galois extension}
\def\ta{\tilde a}
\def\tb{\tilde b}
\def\td{\tilde d}
\def\st{\stackrel}
\def\<{\langle}
\def\>{\rangle}
\def\d{\mbox{$\mathop{\mbox{\rm d}}$}}
\def\id{\mbox{$\mathop{\mbox{\rm id}}$}}
\def\ker{\mbox{$\mathop{\mbox{\rm Ker$\,$}}$}}
\def\hom{\mbox{$\mathop{\mbox{\rm Hom}}$}}
\def\im{\mbox{$\mathop{\mbox{\rm Im}}$}}
\def\map{\mbox{$\mathop{\mbox{\rm Map}}$}}
\def\o{\sp{[1]}}
\def\t{\sp{[2]}}
\def\mo{\sp{[-1]}}
\def\z{\sp{[0]}}

\def\ses{short exact sequence}
\def\csa{$C^*$-algebra}
\def\cO{{\mathcal O}}
\def\O{\cO}
\newcommand{\Boxneu}{\square}
\def\C{{\Bbb C}}
\def\N{{\Bbb N}}
\def\R{{\Bbb R}}
\def\Z{{\Bbb Z}}
\def\T{{\Bbb T}}
\def\cT{{\cal T}}
\def\cK{{\cal K}}
\def\cH{{\cal H}}
\def\im{{\rm Im}}
\def\id{{\rm id}}
\def\tr{{\rm tr}}
\def\Tr{{\rm Tr}}
\def\B{{\mathcal B}}
\def\H{{\mathcal H}}

\title{\large\bf \vspace*{-15mm} A LOCALLY TRIVIAL QUANTUM HOPF FIBRATION}
\author{
\vspace*{-1mm}\large\sc 
Piotr M.~Hajac\\
\vspace*{-2mm}\normalsize 
Mathematisches Institut, Universit\"at M\"unchen\\
\vspace*{-1mm}\normalsize 
Theresienstr.\ 39, M\"unchen, 80333, Germany\\
\vspace*{-1mm}\normalsize 
and\\
\vspace*{-2mm}\normalsize 
Instytut Matematyczny, Polska Akademia Nauk\\
\vspace*{-1mm}\normalsize 
ul.\ \'Sniadeckich 8, Warszawa, 00-950 Poland\\
\vspace*{-1mm}\normalsize 
and\\
\vspace*{-2mm}\normalsize 
Katedra Metod Matematycznych Fizyki, Uniwersytet Warszawski\\
\vspace*{-1mm}\normalsize 
ul.\ Ho\.za 74, Warszawa, 00-682 Poland \vspace{0mm}\\
\large
http://www.fuw.edu.pl/$\!\widetilde{\phantom{m}}\!$pmh
\normalsize 
\and
\vspace*{-1mm}\large\sc
Rainer Matthes\\
\vspace*{-2mm}\normalsize
Fachbereich Physik der TU Clausthal\\
\vspace*{-1mm}\normalsize 
Leibnizstr. 10, D-38678 Clausthal-Zellerfeld, Germany\\
\large
e-mail: ptrm@pt.tu-clausthal.de, matthes@itp.uni-leipzig.de
\and
\vspace*{-1mm}\large\sc 
Wojciech Szyma\'nski\\
\vspace*{-2mm}\normalsize 
School of Mathematical and Physical Sciences, University of Newcastle\\
\vspace*{-1mm}\normalsize 
Callaghan, NSW 2308, Australia\\
\large 
e-mail: wojciech@frey.newcastle.edu.au
}
\date{\normalsize }

\maketitle

\begin{abstract}\normalsize
The irreducible $*$-representations of the polynomial algebra $\cO(S^3_{pq})$ of
the  quantum 3-sphere introduced by Calow and Matthes are classified. 
The $K$-groups of its
universal $C^*$-algebra are shown to coincide with their classical counterparts.
The $U(1)$-action on $\cO(S^3_{pq})$ corresponding for $p=1=q$ to the classical Hopf fibration
is proven to be Galois (free). The thus obtained locally trivial Hopf-Galois
extension is shown to be relatively projective (admitting a strong connection)
and non-cleft. 
The latter is
proven by determining an appropriate Chern-Connes pairing.
\end{abstract}

\section*{Introduction}
\setcounter{equation}{0}

Splitting and gluing topological spaces along 2-spheres or 2-tori are standard
procedures in the study of 3-dimensional manifolds. Fibering such manifolds is
another important tool revealing their geometry. In the case of $S^3$, we have
the well-known Heegaard splittings and Hopf fibration. The former present $S^3$
as two copies of a solid torus glued along their boundaries, and the latter as
a non-trivial principal $U(1)$-bundle over $S^2$.

In \cite{m-k91a}, K.\ Matsumoto applied the idea of a Heegaard splitting
to construct a noncommutative 3-sphere $S^3_\theta$ out of two quantum solid tori.
Then the $U(1)$-action on $S^3_\theta$ was defined and the quotient space
$S^3_\theta/U(1)$ proven to coincide with $S^2$ \cite{m-k91b}. 
Thus a noncommutative
Hopf fibration was constructed. Here we study along these lines a different but 
analogously constructed example of a quantum 
3-sphere \cite{cm}.

Throughout the paper we use the jargon of Noncommutative Geometry referring
to quantum spaces as objects dual  to noncommutative algebras
in the sense of the Gelfand-Naimark correspondence
between spaces and function algebras. 
The unadorned tensor product means the completed tensor product
when placed between $C^*$-algebras (this is not ambiguous as all $C^*$-algebras 
we consider are nuclear), and the algebraic tensor product over
$\C$ otherwise. The algebras are assumed to be associative and over $\C$.
They are also unital unless the contrary is obvious from the context.
$C_0(\mbox{locally compact Hausdorff space})$ means the (non-unital)
algebra of vanishing-at-infinity continuous functions on this space. 
By $\cO(\mbox{quantum space})$ we denote the polynomial algebra of a
quantum space, and by $C(\mbox{quantum space})$ the corresponding $C^*$-algebra. 
By classical points we understand 1-dimensional $*$-representations.
In this paper, the $C^*$-completion ($C^*$-closure) of a $*$-algebra always
means the completion with respect to the supremum norm over all
$*$-representations in bounded operators. 

First, we recall the necessary facts and definitions. This includes the 
construction of 
the  quantum 3-sphere  $S^3_{pq}$ introduced in \cite{cm}, which is also
obtained by gluing two quantum solid tori. Here, however, the noncommutativity
comes from the quantum disc \cite{kl93} rather than the quantum torus (see 
\cite{r-ma90} and references therein). Also, contrary to $S^3_\theta$, the sphere
$S^3_{pq}$ was constructed in the spirit of {\em locally trivial principal bundles}.
Indeed, it can be easily noted that as both the base and the fibre of the Matsumoto
noncommutative Hopf fibration are classical, the noncommutativity of $S^3_\theta$
rules out its local triviality. On the other hand, $S^3_{pq}$ is by construction
a locally trivial quantum $U(1)$-space.

We begin the main part of our paper by classifying the unitary classes of
irreducible $*$-representations of the polynomial algebra $\cO(S^3_{pq})$,
finding a basis of $\cO(S^3_{pq})$, and
defining the $C^*$-algebra $C(S^3_{pq})$.
Then we discuss the fact that $C(S^3_{pq})$ is a graph $C^*$-algebra for
a certain 2-graph.
Next, we prove that 
$K_i(C(S^3_{pq}))\cong\Z\cong K^i(S^3)$, $i\in\{0,1\}$. On the algebraic side,
we show that  $S^3_{pq}$ is a quantum {\em principal} $U(1)$-bundle in the sense
of Hopf-Galois theory. More precisely, by constructing a strong connection, 
we prove that $\cO(S^2_{pq})\inc\cO(S^3_{pq})$ is a
{\em relatively projective} $\cO(U(1))$\gal. From a formula for the strong connection,
we determine the idempotent matrices of all quantum line bundles associated to
the noncommutative Hopf fibration $S^3_{pq}\rightarrow S^2_{pq}$.
 Finally, we pair  a trace on $\cO(S^2_{pq})$
with the idempotent corresponding to the identity representation of $U(1)$.
Since this pairing turns out to be
non-trivial, we conclude that the extension is non-cleft. Thus we have an
example of a locally trivial relatively projective noncommutative \hge\
which not only is not trivial, but also is not a cross-product construction.

Recently, there was an outburst of new examples of noncommutative 3 and 4-spheres. 
For a comparison study of these and older constructions we refer to
\cite{d-l}. Let us only mention that there exist at least two more classes of
non-classical Hopf fibrations, notably 
the quantum Hopf fibrations coming from $SU_q(2)$
and the super Hopf fibration. (See \cite{bm00} and  references
therein for the former, and  \cite{dgh01} and its references for the latter.) 
Although $C(S^3_{pq})$ is not isomorphic to $C(SU_q(2))$ (different
sets of classical points), the \csa\ of the generic Podle\'s quantum sphere
coincides with the \csa\ of the base space of our locally trivial quantum Hopf
fibration, i.e.,  $C(S^2_{pq})\cong C(S^2_{\mu c})$, $p,q,\mu\in(0,1)$, $c>0$ 
\cite[Proposition~21]{cm00}. 
This holds despite the fact that the polynomial $*$-subalgebras $\O(S^2_{pq})$
and  $\O(S^2_{\mu c})$ are non-isomorphic, as shown in Proposition~\ref{non}.
\section{Preliminaries}
\setcounter{equation}{0}

\subsection{Locally trivial \boldmath$H$-extensions}
\label{lotri}

The idea of a locally trivial $H$-extension can be traced back to
\cite{bm93,p-mj94,d-m96,bk96,cm}. The prerequisite idea of the gluing of
quantum spaces can be traced much further. In terms of \csa s, the gluing 
corresponds to the pullback construction, which is essential, e.g., for
the Busby invariant. We refer to \cite{p-gk99} for lots of generalities
on pullbacks of \csa s. Here, let us only recall the needed definitions
and fix the terminology.
\bde(\cite{bk96},\cite{cm00})
A covering of an algebra $B$ is a family $\{J_i\}_{i\in I}$ of ideals with zero
intersection. Let $\pi_i:B\ra B_i:=B/J_i$, $\pi^i_j:B_i\ra B_{ij}:=
B/(J_i+J_j)$ be the quotient maps.
A covering $\{J_i\}_{i\in I}$ is called {\em complete} iff the homomorphism
$$
B\ni b\longmapsto (\pi_i(b))_{i\in I}\in B_c:=
\left\{(b_i)_{i\in I}\in\mbox{$\prod_{i\in I}$}B_i\;
|\;\pi^i_j(b_i)=\pi^j_i(b_j)\right\}
$$
is surjective. (It is automatically injective.) 
\ede
Note that finite coverings consisting of closed ideals in $C^*$-algebras and two-element
coverings are always complete.
See \cite{cm00} for more information
about completeness of coverings, in particular for an example of a
non-complete covering.
Let $B$ and $P$ be algebras and $H$ be a Hopf algebra. 
Assume that $B$ is a subalgebra of $P$
via an injective homomorphism $\iota:B\ra P$, which we also write as $B\inc P$.
The inclusion 
$B\inc P$ is called an $H$-extension if $P$ is a right $H$-comodule algebra via
a right coaction $\Delta_R:P\ra P\ot H$ such that $B$ coincides with the subalgebra
of coinvariants, i.e.,  $B=P^{coH}:=\{p\in P\;|\;\Delta_R(p)=p\ot 1\}$.
\bde\label{lth}
(\cite{cm})
A {\em locally trivial} $H$-extension\footnote{
In \cite{cm} locally trivial $H$-extensions were called
locally trivial quantum principal fibre bundles.
}
 is an $H$-extension $B\inc P$ together
with the following (local) data:\\[.3cm]
(i) a complete finite covering $\{J_i\}_{i\in I}$ of $B$;\\[.3cm]
(ii) surjective homomorphisms
$\chi_i:P\ra B_i\ot H$ (local trivializations) such that\\[.3cm]
\phantom{(ii)} (a) $\chi_i\ci\iota=\pi_i\ot 1$,\\[.3cm]
\phantom{(ii)} (b) $(\chi_i\ot\id)\ci\Delta_R=(\id\ot\Delta)\ci\chi_i$,\\[.3cm]
\phantom{(ii)} (c) $\{\ker\chi_i\}_{i\in I}$ is a complete covering of $P$.
\ede
Locally trivial $H$-extensions can be reconstructed from the local data.
First, it follows from the proof of \cite[Proposition 2]{cm} that there
exist homomorphisms $\phi_{ij}:B_{ij}\ot H\ra B_{ij}\ot H$ (to be thought
of as the change of local trivializations) that are uniquely determined 
by the formula
\[
\varphi_{ij}\ci(\pi^j_i\ot\id)\ci\chi_j=(\pi^i_j\ot\id)\ci\chi_i.
\]
As argued in \cite{cm} (see formulas (5) and (6) and the remark after the proof of 
Proposition~2),
these maps
are isomorphisms satisfying 
\[
(\id\ot\Delta)\ci\varphi_{ij}=(\varphi_{ij}\ot\id)\ci(\id\ot\Delta),
\]
\[
\varphi_{ij}(b\ot 1)=(b\ot 1).
\]
Next, in analogy with the
classical situation, one defines transition functions
$\tau_{ji}:H\ra B_{ij}$ by 
\[
\tau_{ji}(h)=(\id\ot\varepsilon)(\varphi_{ij}(1\ot h)).
\]
Equivalently, one has
\beq\label{phi}
\phi_{ij}(b\ot h)=b\tau_{ji}(h_{(1)})\ot h_{(2)}.
\eeq
It can be shown (see \cite{bk96} or \cite[Proposition 4]{cm}) that the transition functions of a 
locally trivial $H$-extension are homomorphisms
with the following properties:
\bea
\label{tii} \tau_{ii}&=&\varepsilon,\\
 \tau_{ji}\ci S&=&\tau_{ij},\\
 \tau_{ij}(H)&\subset &Z(B_{ij})\;\;\;\mbox{(centre of $B_{ij}$)},\\
\label{tijk} \pi^{ij}_k\ci\tau_{ij}&=&m_{B_{ijk}}\ci((\pi^{ik}_j\ci\tau_{ik})
\ot(\pi^{jk}_i\ci\tau_{kj}))\ci\Delta.
\eea
Here $\pi^{ij}_k:B_{ij}\ra B_{ijk}:=B/(J_i+J_j+J_k)$ are the quotient maps
and $m_{B_{ijk}}$ are the
multiplications in $B_{ijk}$. 

Conversely, let us consider an algebra $B$ with a complete finite covering
$\{J_i\}_{i\in I}$ and a Hopf algebra $H$. Assume that we have a family of
homomorphisms $\tau_{ji}:H\ra B_{ij}$ satisfying (\ref{tii})--(\ref{tijk}).
Define $\phi_{ij}$ by the formula (\ref{phi}) and put
\beq\label{pglu}
\breve{P}
=
\left\{(f_i)_{i\in I}\in\mbox{$\prod_{i\in I}$}(B_i\ot H)\;
|\;(\pi^i_j\ot\id)(f_i)=\varphi_{ij}((\pi^j_i\ot\id)(f_j))\right\}.
\eeq
One can verify that the formulas 
\bea\label{drglu}
\breve{\Delta}_R((f_i)_{i\in I})&=&((\id\ot\Delta)(f_i))_{i\in I},
\\
\breve{\chi}_i((f_i)_{i\in I})&=&f_i,
\\
\label{iglu}
\breve{\iota}(b)&=&(\pi_i(b)\ot 1)_{i\in I}\,,
\eea
turn $\breve P$ into a locally trivial $H$-extension of $B$.
Moreover, we have:
\bpr \label{bueglu}(\cite{cm})
Let $P$ be a locally trivial $H$-extension of $B$ corresponding to covering 
$\{J_i\}_{i\in I}$, and  $\tau_{ji}:H\ra B_{ij}$ be its transition functions.
Let $\breve P$ be a locally trivial $H$-extension constructed from $\tau_{ji}$'s.
Then the formula $p\mapsto(\chi_i(p))_{i\in I}$ defines an isomorphism of
locally trivial $H$-extensions $P$ and $\breve P$.
\epr

\subsection{Construction of \boldmath$S^3_{pq}$}

Our starting point is 
the coordinate algebra ${\cal O}(D_p)$ of the quantum disc, which is defined as the
universal unital $*$-algebra generated by  $x$ fulfilling the relation
\beq \label{redi}
x^*x-pxx^*=1-p,\;\;\;0<p<1.
\eeq
This is a one-parameter sub-family of the two-parameter family 
of quantum discs defined in \cite{kl93}. Using all bounded representations of
 ${\cal O}(D_p)$,
one can define its $C^*$-closure $C(D_p)$. It can be shown that the 
\csa\ $C(D_p)$
is isomorphic with the Toeplitz algebra $\cT$ (e.g., 
see \cite[Proposition 15]{cm00}), and that $\|x\|=1$ (see 
\cite[Proposition IV.1(I)]{kl93}). Let us also mention that there are  
unbounded representations of the relation (\ref{redi}). 
They are given, e.g., in \cite[Section 5.2.6]{ks97}.

Next, we glue two quantum discs to get a quantum $S^2$. Let $\cO(S^1)$ be the
universal $*$-algebra generated by the unitary $u$. Then we have a natural
epimorphism $\pi_p:{\cal O}(D_p)\ra \cO(S^1)$ given by
$\pi_p(x)=u$. (This corresponds to embedding $S^1$ into $D_p$ as its boundary.)
Now, we can define 
$\cO(S^2_{pq})$ in the following manner \cite{cm00}:
\beq\label{s2glu}
\cO(S^2_{pq}):=\{(f,g)\in \cO(D_p)\oplus\cO(D_q)~|~
\pi_p(f)=\pi_q(g)\}.
\eeq
$\cO(S^2_{pq})$ has a complete covering $\{\ker pr_1,\ker pr_2\}$, where
$pr_1$ and $pr_2$ are the restrictions to $\cO(S^2_{pq})$ of the  projections
on ${\cal O}(D_p)$ and ${\cal O}(D_q)$, respectively
(see \cite[Proposition 8]{cm00}). Furthermore,
one has canonical isomorphisms $\cO(S^2_{pq})/\ker pr_1\cong \cO(D_q)$,
$\cO(S^2_{pq})/\ker pr_1\cong \cO(D_q)$ and 
$\cO(S^2_{pq})/(\ker pr_1+\ker pr_2)\cong \cO(S^1)$.
As was shown in \cite[Proposition 17]{cm00}, 
$\cO(S^2_{pq})$ can be identified with the universal
$*$-algebra generated by $f_0$ and $f_1$ satisfying the relations
\bea
\label{ref1} f_0&=&f_0^*,\\
f_1^*f_1-qf_1f_1^*&=&(p-q)f_0+1-p,\\
f_0f_1-pf_1f_0&=&(1-p)f_1,\\
\label{ress}(1-f_0)(f_1f_1^*-f_0)&=&0.
\eea
The isomorphism is given by $f_1\mapsto (x,y)$, $f_0\to (xx^*,1)$. 
Here $x$ denotes
the generator of the $*$-algebra $\cO(D_p)$, and $y$ that of $\cO(D_q)$.
In terms of these generators, the irreducible $*$-representations of
$\cO(S^2_{pq})$ can be given as follows \cite[Proposition 19]{cm00}:
\bea
&&
\rho_\theta(f_0)=1,\;\;\;\rho_\theta(f_1)=e^{i\theta},\;\;\;\theta\in[0,2\pi) 
\;\;\;\mbox{(classical points),}
\\ &&\label{ro1}
\rho_1(f_0)e_k=(1-p^k)e_k,\;\;\;
\rho_1(f_1)e_k=\sqrt{1-p^{k+1}}\:e_{k+1},\;\;\;k\geq 0;
\\ &&\label{ro2}
\rho_2(f_0)e_k=e_k,\;\;\;
\rho_2(f_1)e_k=\sqrt{1-q^{k+1}}\:e_{k+1},\;\;\;k\geq 0.
\eea
Here $\{e_k\}_{k\geq 0}$ is an orthonormal basis of a separable Hilbert space.

In the classical case $p=q=1$, the relations (\ref{ref1})--(\ref{ress}) reduce
to commutativity and the geometrical relation (\ref{ress}). Adding by hand
the conditions $|f_0|\leq 1,\ |f_1|\leq 1$ (which are automatic in the
noncommutative case \cite[Proposition 19]{cm00}), 
one obtains as the corresponding geometric space a closed
 cone. The irreducible representations given above allow one
to make an analogous picture also in the noncommutative case.
The sum of the squares of the hermitian generators
$f_+=\frac{1}{2}(f_1+f_1^*)$ and $f_-=\frac{i}{2}(f_1-f_1^*)$ is diagonal
in the representations, and one can imagine a discretized version
of the above cone, with the edge being the circle of classical points
and the remainder of the cone being formed by ``non-classical circles'', 
accumulating at this edge (cf.\ \cite[pp.337-8]{cm00}).

It follows from the relations (\ref{ref1}) -- (\ref{ress}) that
$\|\rho(f_i)\|\leq 1$ in every bounded representation of $\cO(S^2_{pq})$.
Therefore, one can form a $C^*$-algebra $C(S^2_{pq})$ using bounded
representations.
As noticed already at the end of the introduction, this $C^*$-algebra is 
isomorphic to the $C^*$-algebra $C(S^2_{\mu c})$ of the Podle\'s spheres
for $c>0$. Such an isomorphism does not exist on the level of polynomial
$*$-algebras. Indeed, there are two infinite-dimensional irreducible 
representations $\pi_+$
and $\pi_-$ of $C(S^2_{\mu c})$ \cite[Proposition 4]{p-p87} whose 
restrictions to $\cO(S^2_{\mu c})$
are faithful. (Their faithfulness can be proved by direct arguments using a vector
space basis of $\cO(S^2_{\mu c})$ provided in
\cite[Section 3]{p-p89}).) On the other hand, the $*$-algebra 
$\cO(S^2_{pq})$ has only
two inequivalent infinite-dimensional irreducible representations 
$\rho_1$ and $\rho_2$ \cite{cm00}. None of them is faithful, since $1-f_0$ is in 
the kernel of $\rho_2$ and $f_1f_1^*-f_0$
is in the kernel of $\rho_1$. Thus, we have
\bpr\label{non}
For any $p,q,\mu\in (0,1),~c>0$, there is no $*$-algebra isomorphism
between $\cO(S^2_{pq})$ and $\cO(S^2_{\mu c})$.
\epr

We are now ready for the definition of $\cO(S^3_{pq})$.
We consider $\cO(U(1))$ as  the $*$-algebra $\cO(S^1)$ equipped with the
 Hopf algebra structure
given by $\Delta(u)=u\ot u$, $\varepsilon(u)=1$, $S(u)=u^*$,
and $\cO(S^2_{pq})$ as an algebra with the complete covering $\{\ker pr_1,\ker pr_2\}$.
The homomorphisms $\tau_{12}:=\id:\cO(U(1))\ra\cO(S^1)$, $\tau_{21}:=S$,
$\tau_{11}:=\he=:\tau_{22}$, evidently fulfill the axioms (\ref{tii})--(\ref{tijk}).
Therefore, we can proceed along the 
lines of Subsection~\ref{lotri}, and define
the following locally trivial $\cO(U(1))$-extension:
\bde(\cite[p.152]{cm})
Let $\tau_{ji}:\cO(U(1))\ra\cO(S^2_{pq})_{ij}$ be the homomorphisms given above.
We define $\cO(S^3_{pq})$ as the locally trivial $\cO(U(1))$-extension of 
$\cO(S^2_{pq})$ given by $\tau_{ji}$'s via (\ref{phi}) and (\ref{pglu}). 
Explicitly, the algebra $\cO(S^3_{pq})$ is 
$$
\left\{(a_1,a_2)\in 
(\cO(D_p)\ot\cO(U(1))\oplus (\cO(D_q)\ot\cO(U(1))\;|\;
(\pi_p\ot id)(a_1)=\phi_{12}((\pi_q\ot id)(a_2))\right\},
$$
where $\pi_p:\cO(D_p)\ra \cO(S^1)$, $\pi_p(x)=u$, 
$\pi_q:\cO(D_q)\ra \cO(S^1)$, $\pi_q(y)=u$. 
\ede
Note that we glue two  quantum solid tori $D_p\times U(1)$ and $D_q\times U(1)$
along their classical boundaries, which are $T^2$. 
The subspace of the classical points of the resulting $S^3_{pq}$ 
is precisely the locus of
the gluing (see Section~2). In terms of generators and relations,
$\cO(S^3_{pq})$ can be characterised 
 in the following way:
\ble(\cite[Proposition 21]{cm})
\label{s3}
$\cO(S^3_{pq})$ is isomorphic to the universal 
unital $*$-algebra generated by $a$ and $b$
satisfying the relations
\begin{eqnarray}
\label{prel1}{a}^*{a} - q{a} {a}^* &=& 1-q, \\
\label{prel2} {b}^* {b}-p{b} {b}^* &=& 1-p, \\
\label{prel3} ab=ba,~~ a^*b&=&ba^*, \\
\label{prel4} (1-{a}{a}^*)(1-{b} {b}^*)&=& 0. 
\end{eqnarray}
\ele
The isomorphism is given by $(1\ot u,y\ot u)\mapsto a$ and 
$(x\ot u^*,1\ot u^*)\mapsto b$. 
Using this identification and the description of $\cO(S^2_{pq})$
in terms of the generators $f_0$ and $f_1$ ((\ref{ref1})--(\ref{ress})),
we can write the structural $*$-homomorphisms of the
locally trivial $\cO(U(1))$-extension $\cO(S^2_{pq})\inc\cO(S^3_{pq})$
in the form
\[
 \label{dpa}\Delta_R(a) = a \otimes u,
\;\;\;\Delta_R(b) = b \otimes u^*,
\]
\[
\chi_p(a)= 1 \otimes u,\;\;\;\chi_q(a) = y \otimes u, \;\;\;
\chi_p(b)= x \otimes u^*,\;\;\;\chi_q(b) = 1 \otimes u^*, 
\]
\[
\label{coin}\iota(f_1) = ba,\;\;\;\iota(f_0)=bb^*.
\]

As shown in Section~\ref{irrep}, we can define the $C^*$-algebra of
$\cO(S^3_{pq})$ as the universal unital $C^*$-algebra generated by $a$ and $b$
satisfying the relations (\ref{prel1})--(\ref{prel4}) and $\| a\|=1=\|b\|$.
The last condition follows from the remaining ones for $p,q\in (0,1)$
for the same reason it is automatically true for the generator of $\cO(D_p)$
(see \cite[Proposition IV.1(I)]{kl93}), but for 
$p=1=q$ we need to put it by hand. In the classical case ($p=1=q$), 
this $C^*$-algebra
coincides with $C(S^3)$ \cite[p.334]{m-k91a}. By the universality of 
the $C^*$-algebra
$C(S^3_{pq})$, the right coaction $\dr$ extends to $C(S^3_{pq})$ and is 
equivalent to a $U(1)$-action on $C(S^3_{pq})$ (see 
\cite[Proposition T.5.21]{w-ne93}). The latter
 reduces for $p=1=q$ to the $U(1)$-action on $S^3$ yielding
the Hopf fibration. (Our convention for the action differs from the convention in 
\cite{m-k91a}.)
To understand more precisely the classical case (see \cite[Section~0.3]{n-gl97} 
for related details), 
let us prove the following:
\bpr\label{classic}
Define $X=\{(z_1,z_2)\in\C^2\;|\;(1-|z_1|^2)(1-|z_2|^2)=0,\; |z_i|\leq 1\}$
and $S^3=\{(c_1,c_2)\in\C^2\;|\; |c_1|^2+|c_2|^2=1\}$. The group $U(1)$ acts
on $X$ and $S^3$ via $(z_1,z_2)\cdot e^{i\phi}=(z_1e^{i\phi},z_2e^{-i\phi})$
and $(c_1,c_2)\cdot e^{i\phi}=(c_1e^{i\phi},c_2e^{i\phi})$, respectively, and $X$
and $S^3$ are homeomorphic as $U(1)$-spaces.
\epr
\bpf
The $U(1)$-action on $\C^2$ clearly restricts to both $X$ and $S^3$. Note first
that we can equivalently write the equation $(1-|z_1|^2)(1-|z_2|^2)=0$ in the
form $|z_1|^2+|z_2|^2=1+|z_1|^2|z_2|^2$. This suggests that we can define a map
$X\st{f}{\ra}S^3$ by the formula (cf.\ \cite[p.38]{mt92} 
for the case of $S^3_\theta$):
\[\label{f}
f((z_1,z_2))=(|z_1|^2+|z_2|^2)^{-\frac{1}{2}}(z_1,\overline{z_2}).
\]
Indeed, $f$ is a continuous $U(1)$-map into $S^3$. To find the inverse of $f$, we
look for a map of the form $(c_1,c_2)\to\ha(c_1,\overline{c_2})$.\footnote{
We are grateful to Andrzej Sitarz for his involvement here.
}
A direct computation provides us with the formula:
\[\label{g}
g((c_1,c_2))=\frac{\sqrt{2}(c_1,\overline{c_2})}{\sqrt{1+|\;2|c_1|^2-1\;|}}
=:(g_1,g_2).
\]
Taking advantage of $|c_1|^2+|c_2|^2=1$, we compute
\[\label{gg}
|g_1|^2+|g_2|^2
=
\frac{2|c_1|^2}{1+|\;2|c_1|^2-1\;|}+\frac{2|c_2|^2}{1+|\;2|c_1|^2-1\;|}
=
\frac{2}{1+|\;2|c_1|^2-1\;|}
\]
and
\[
1+|g_1|^2|g_2|^2
=
\frac{(1+|\;2|c_1|^2-1\;|)^2+4|c_1|^2(1-|c_1|^2)}{(1+|\;2|c_1|^2-1\;|)^2}
=
\frac{2}{1+|\;2|c_1|^2-1\;|}.
\]
Hence $(1-|g_1|^2)(1-|g_2|^2)=0$. Furthermore, as 
$|\;2|c_1|^2-1\;|=|\;2|c_2|^2-1\;|$,
\[
|g_i|\leq 1\Leftrightarrow
2|c_i|^2\leq 1+|\;2|c_i|^2-1\;|.
\]
 We have two cases. For $2|c_i|^2\geq 1$, the latter
inequality reads $0\leq 0$. Otherwise, i.e., for $2|c_i|^2< 1$, it is the same as
$2|c_i|^2\leq1$. Thus we have a continuous map $S^3\st{g}{\ra}X$. As both $f$ and
$g$ are evidently $U(1)$-maps, it only remains to prove that they are mutually 
inverse. Remembering (\ref{gg}), we have:
\[
(f\ci g)((c_1,c_2))
=
\frac{\sqrt{2}(c_1,c_2)}{\sqrt{1+|\;2|c_1|^2-1\;|}}
\frac{\sqrt{1+|\;2|c_1|^2-1\;|}}{\sqrt{2}}
=
(c_1,c_2).
\] 
For the other identity, note first that, due to $(1-|z_1|^2)(1-|z_2|^2)=0$,
we have $|z_j|^2\leq|z_i|^2=1$.
To avoid confusion, let us fix $|z_2|^2\leq|z_1|^2=1$.
 (The other case behaves in the
same way.) Now we can compute:
\bea
&&\!\!\!\!\!\!\!\!\!(g\ci f)((z_1,z_2))
\\ &=&
\frac{\sqrt{2}}{\sqrt{|z_1|^2+|z_2|^2}}
\left(
\frac{z_1}{\sqrt{1+\left|\frac{2|z_1|^2}{|z_1|^2+|z_2|^2}-1\right|}},
\frac{z_2}{\sqrt{1+\left|\frac{2|z_2|^2}{|z_1|^2+|z_2|^2}-1\right|}}
\right)
\\ &=&
\left(
\frac{\sqrt{2}z_1}{\sqrt{|z_1|^2+|z_2|^2+|\;2|z_1|^2-|z_1|^2-|z_2|^2\;|}},
\frac{\sqrt{2}z_2}{\sqrt{|z_1|^2+|z_2|^2+|\;2|z_2|^2-|z_1|^2-|z_2|^2\;|}}
\right)\phantom{mmm}
\\ &=&
\frac{\sqrt{2}(z_1,z_2)}{\sqrt{|z_1|^2+|z_2|^2+|z_1|^2-|z_2|^2}}
\\ &=&
(z_1,z_2).
\eea
This ends the proof.
\epf

\subsection{Hopf-Galois extensions and associated modules}

We refer to \cite{m-s93} for  generalities concerning Hopf-Galois theory.
Recall that an $H$-extension $B\inc P$
is called Hopf-Galois iff the canonical map 
\[
can: P\ot_BP\lra P\ot H,\;\;\;p\ot p'\mapsto p\dr(p'),
\]
is bijective. Note that $can$ is surjective whenever, for any generator $h$ of $H$,
the element $1\ot h$ is in the image of $can$ (cf.\ \cite[pp.106--7]{s-p00}). Indeed, if 
$\sum_ih_i\ot\tilde{h}_i,\,\sum_jg_j\ot\tilde{g}_j\in P\ot P$ are tensors
such that $can(\sum_ih_i\ot_B\tilde{h}_i)=1\ot h$ and 
$can(\sum_jg_j\ot_B\tilde{g}_j)=1\ot g$, then 
$\sum_{ij}g_jh_i\ot\tilde{h}_i\tilde{g}_j\in P\ot P$ has the property
$can(\sum_{ij}g_jh_i\ot_B\tilde{h}_i\tilde{g}_j)=1\ot hg$.
Hence, for any monomial $w\in H$, the element $1\ot w$ is in the image of $can$,
and its surjectivity follows from its left $P$-linearity.

 The restricted inverse of $can$, $T:=can^{-1}\ci(1\ot\id)$,
is called the translation map. We are interested in unital bicolinear liftings
of the translation map $T$ because, when the antipode $S$ of $H$ is bijective,
they can be interpreted as strong connections on algebraic quantum principal bundles
\cite{bh}. More precisely, we take the canonical surjection $\pi_B$ 
and demand that the
following diagram be commutative: 
\begin{diagram}[height=10mm]
  &             & P\ot P \\
  & \ruTo{\ell} & \dTo_{\pi_B} \\
H & \rTo{T}     & P\ot_BP. 
\end{diagram}
Then we equip $P\ot P$ with an $H$-bicomodule structure via the maps
\[
\hD_L^\otimes:=\llp (S^{-1}\ot\id)\ci(\mbox{flip})\ci\dr\lrp\ot\id\;\;\;
\mbox{and}\;\;\;\hD_R^\otimes:=\id\ot\dr,
\]
and require that
\[
\hD_L^\otimes\ci\ell=(\id\ot\ell)\ci\hD\;\;\;\mbox{and}\;\;\;
\hD_R^\otimes\ci\ell=(\ell\ot\id)\ci\hD.
\]
Finally, we ask that $\ell(1)=1\ot 1$. In the more general context of symmetric
coalgebra-Galois extensions, it is shown in \cite{bh} that the existence of such
a lifting is equivalent to $P$ being an $(H^*)^{op}$-relatively projective left 
$B$-module. (Here $(H^*)^{op}$ means the convolution algebra of functionals on $H$
taken with the opposite multiplication.) We call such \hge s relatively projective.
In general (cf.\ 
\cite[p.197]{ce56}), we say
that a $(B,A)$-bimodule $P$ is an {\em $A$-relatively projective left $B$-module} iff
   for every diagram with the exact row
\begin{diagram}[height=10mm]
M & \pile{\rTo^{\pi}\\ \lTo_{i}} &  N       & \rTo & 0 \\
  &                              & \uTo_{f} &      &    \\
  &                              & P,       &      &
\end{diagram}
where $M,N$ are $(B,A)$-bimodules, $\pi$ and $f$ are $(B,A)$-bimodule maps
and $i$ is a right $A$-module splitting of $\pi$, there exists a
$(B,A)$-bimodule map $g$ rendering the following diagram commutative
\begin{diagram}[height=10mm]
M & \pile{\rTo^{\pi}\\ \lTo_{i}} &  N       & \rTo & 0 \\
  &           \luTo_g            & \uTo_{f} &      &    \\
  &                              & P.       &      &
\end{diagram}

If $B\inc P$ is an Hopf-Galois $H$-extension and $\rho:V\ra V\ot H$ is a coaction, then we can define
the associated left $B$-module $\hom_\rho(V,P)$ of all colinear maps from $V$ to $P$.
Such modules play the geometric role of the modules of sections of associated vector
bundles. Therefore, to be in line with the Serre-Swan theorem and $K$-theory, it is
desirable to have them finitely generated projective. It turns out that this is 
always the case for $\dim V<\infty$ and relatively projective \hge s with
bijective antipodes (see \cite[Corollary~2.6]{dgh01}, cf.\ \cite{bh} for a more 
general context).

Assume that $\rho$ is a 1-dimensional corepresentation. Then it is given by a 
group-like $g$, $\rho(1)=1\ot g$. Assume further that $B\inc P$ is an $H$\gal\
admitting a strong connection $\ell$ (relative projectivity).
Put $l(g)=\sum_{k=1}^nl_k(g)\ot r_k(g)$.
Then it can be shown \cite{bh} that $E_{jk}:=r_j(g)l_k(g)\in B$, 
the matrix $E:=(E_{jk})$ is idempotent
($E^2=E$), and we have an isomorphism of left 
$B$-modules $\hom_\rho(\C,P)\cong B^nE$.

\section{Representations and the \boldmath$C^*$-algebra of
\boldmath $\cO(S^3_{pq})$}
\setcounter{equation}{0}

To begin with, we classify the bounded irreducible
$*$-representations of $\cO(S^3_{pq})$. We do it much as in the case of the
quantum real projective space $\R P^2_q$\cite[Theorem 4.5]{hms}. 
As a result, we obtain two
$S^1$-families of infinite-dimensional representations and a $T^2$-family
of one-dimensional representations (classical points). The latter proves
that  $S^3_{pq}$ differs from other quantum 3-spheres (see \cite{d-l} for
details).
\bth\label{irrep}
Let $\cH$ be a separable Hilbert space with an orthonormal basis
$\{e_k\}_{k\geq 0}$. Any irreducible $*$-representation of $\cO(S^3_{pq})$
in bounded operators
on a Hilbert space
is unitarily equivalent to one of the following:
\[\label{r1}
\rho_{1\theta}(a)e_k=e^{i\theta}e_k,\;\;\;
\rho_{1\theta}(b)e_k=\sqrt{1-p^{k+1}}\:e_{k+1},\;\;\;
\theta\in[0,2\pi);
\]
\[\label{r2}
\rho_{2\theta}(a)e_k=\sqrt{1-q^{k+1}}\:e_{k+1},\;\;\;
\rho_{2\theta}(b)e_k=e^{i\theta}e_k,\;\;\;
\theta\in[0,2\pi);
\]
\[\label{r3}
\rho_{\theta_1\theta_2}(a)=e^{i\theta_1},\;\;\;
\rho_{\theta_1\theta_2}(b)=e^{i\theta_2},\;\;\;
\theta_1,\theta_2\in[0,2\pi).
\]
\ethe
\bpf
Let $\rho$ be a $*$-representation of $\cO(S^3_{pq})$ in bounded operators 
on a Hilbert space $\widetilde\cH$. 
Then it follows immediately from the relations (\ref{prel1})--(\ref{prel4}) that 
$\ker(1-\rho(aa^*))$ and $\ker(1-\rho(bb^*))$ are invariant subspaces.
Hence $\cH_0:= \ker(1-\rho(aa^*))\cap \ker(1-\rho(bb^*))$ is invariant.
Let  $\phi\neq 0$ be a vector in the orthogonal complement of the closure 
of the sum $\ker(1-\rho(aa^*))+\ker(1-\rho(bb^*))$. Then, due to the
invariance of this complement, $(1-\rho(aa^*))(1-\rho(bb^*))\phi\neq 0$,
which contradicts (\ref{prel4}). Therefore this complement must be zero,
i.e., $\widetilde\cH$ automatically coincides with the closure of
$\ker(1-\rho(aa^*))+\ker(1-\rho(bb^*))$. Thus we have the direct sum
decomposition 
\beq
\widetilde\cH=\cH'\oplus \cH''\oplus \cH_0
\eeq
into $\rho$-invariant subspaces. Here $\cH':=\ker(1-\rho(aa^*))\ominus \cH_0$ and
$\cH'':=\ker(1-\rho(bb^*))\ominus \cH_0$ are appropriate orthogonal complements.
Our strategy is to look for irreducible representations on these
three subspaces separately. We abuse notation by using $\rho$ to denote also
its restrictions.

For the restriction of $\rho$ to $\cH_0$ we have 
$\rho(a)\rho(a^*)=1=\rho(b)\rho(b^*)$. 
It follows from the disc-like relations (\ref{prel1})
and (\ref{prel2}) that also $\rho(a^*)\rho(a)=1=\rho(b^*)\rho(b)$, so that
$\rho(a)$ and $\rho(b)$ are unitary. Since $\rho(a)$ and $\rho(b)$
commute, we arrive at the third family of representations.

Consider now the restriction of $\rho$ to $\cH'$. On this invariant subspace, we
have $\rho(aa^*)=1$, and it follows from (\ref{prel1}) that $\rho(a)$ is
unitary. The operator $1-\rho(bb^*)$ is injective on $\cal H'$.
Our aim is to determine the spectrum of $1-\rho(bb^*)\in B(\cH')$.
First we show that
$spec(1-\rho(bb^*))\subset [0,1]$. From (\ref{prel2}) we can conclude
that $\|\rho(b^*b-pbb^*)\|=1-p$.
Hence $\|\rho(b^*b)\|-p\|\rho(bb^*)\|\leq1-p$, which gives
\beq
\|\rho(b)\|^2=\|\rho(b^*b)\|=\|\rho(bb^*)\|\leq 1.
\eeq
This means $0\leq \rho(bb^*)\leq 1$, and therefore also 
$0\leq 1-\rho(bb^*)\leq 1$, which yields the desired inclusion for the 
spectrum. 
If $0\in spec(1-\rho(bb^*))$, it cannot be an eigenvalue, since this would contradict
the injectivity of $1-\rho(bb^*)$. Therefore, 0 cannot be isolated in the spectrum. 
If $1$ were the only element of the spectrum, this
would mean $1-\rho(bb^*)=1$ on $\cH'$, i.e., $\rho(b)=0$, contradicting 
(\ref{prel2}). Hence we conclude that there exists 
$\lambda\in (0,1)\cap spec(1-\rho(bb^*))$. By \cite[Lemma 3.2.13]{kr97}, 
there exists a sequence of unit vectors
$\{\phi_n\}_{n\in\N}\subset \cH'$ such that 
\beq\label{aeig}
\lim_{n\to\infty}\|\rho(1-bb^*)\phi_n - \lambda\phi_n\|=0.
\eeq
On the other hand, we have the estimate
\bea\label{est}
\|\rho((1-bb^*)b^*)\phi_n-p^{-1}\lambda\rho(b^*)\phi_n\|&=&
\|p^{-1}\rho(b^*(1-bb^*))\phi_n - p^{-1}\lambda\rho(b^*)\phi_n\|\nonumber\\
&\leq& p^{-1}\|\rho(b^*)\|\|\rho(1-bb^*)\phi_n - \lambda\phi_n\|.
\eea
To show that $\|\rho(b^*)\phi_n\|\geq C$ for some $C>0$ and 
all sufficiently big $n$, we compute
\beq
\|\rho(bb^*)\phi_n\|=\|(1-\lambda-\rho(1-bb^*-\lambda))\phi_n\|
\geq |1-\lambda|-\|\rho(1-bb^*-\lambda)\phi_n\|.
\eeq
Using $\|\rho(b)\|\|\rho(b^*)\phi_n\|\geq\|\rho(bb^*)\phi_n\|$,
we  conclude that
\beq
\|\rho(b^*)\phi_n\|\geq\frac{1}{\|\rho(b)\|}
(|1-\lambda|-\|\rho(1-bb^*-\lambda)\phi_n\|).
\eeq
Due to (\ref{aeig}), the term in parentheses  approaches $1-\lambda>0$ for 
$n\ra\infty$, so that there is $N\in\N$ and $C>0$
 such that $\|\rho(b^*)\phi_n\|\geq C$ for $n>N$. 
Consequently,
 we can form a sequence $\{\eta_n\}_{n\in\N}$ of unit vectors consisting of 
$\eta_n:=\frac{\rho(b^*)\phi_n}{\|\rho(b^*)\phi_n\|}$ for $n>N$ and arbitrarily
chosen unit vectors $\eta_n$ for $n\leq N$. It is now immediate from 
(\ref{aeig}) and the estimate (\ref{est}) that
\beq
\lim_{n\to\infty}\|\rho(1-bb^*)\eta_n-p^{-1}\lambda\eta_n\|=0.
\eeq
Employing again \cite[Lemma 3.2.13]{kr97}, we have
 that $p^{-1}\lambda\in spec(\rho(1-bb^*))$. 
We can iterate this reasoning until $p^{-k}\hl=1$. This has to be true for some
$k$, as otherwise we would contradict $spec(\rho(1-bb^*))\inc[0,1]$.
 It follows
that 
\[
\{1,p,\ldots,p^k\}\subset spec(\rho(1-bb^*))\inc
\{1,p,p^2,\ldots\}\cup\{0\}.
\]
We now show that the latter inclusion is an equality.
Let $\xi_k$ be a (non-zero) eigenvector corresponding to the eigenvalue $p^k$, 
$\rho(1-bb^*)\xi^k=p^k\xi_k$. Then, using (\ref{prel2}), we have 
\[
\rho(1-bb^*)\rho(b)\xi_k=p\rho(b)\rho(1-bb^*)\xi_k=p^{k+1}\rho(b)\xi_k.
\]
Using the same relation $(1-bb^*)b=pb(1-bb^*)$, we obtain 
\bea
\|\rho(b)\xi_k\|^2
&=&
\langle \rho(b)\xi_k\,|\,\rho(b)\xi_k\rangle
\nonumber\\ &=&
\langle \rho(b^*b)\xi_k\,|\,\xi_k\rangle
\nonumber\\ &=&
\langle(p\rho(bb^*)+1-p)\xi_k\,|\,\xi_k\rangle
\nonumber\\ &=&
p\|\rho(b^*)\xi_k\|^2+(1-p)\|\xi_k\|^2>0.
\eea
Hence $\rho(b)\xi_k$ is a non-zero eigenvector to the eigenvalue $p^{k+1}$.
This proves that 
\[
spec(\rho(1-bb^*))=\{1,p,p^2,\ldots\}\cup\{0\}.
\]

We are ready now to construct a set of orthonormal vectors.
Since 1 is isolated in $spec(\rho(1-bb^*))$, there exists a normalised eigenvector
$\xi$ given by  $\rho(1-bb^*)\xi=\xi$.
Making again the use of $(1-bb^*)b=pb(1-bb^*)$, we obtain
\bea
\|\rho(b^{k+1})\xi\|^2
&=& 
\langle \rho(b^*bb^k)\xi\,|\,\rho(b^k)\xi\>
\nonumber\\ &=&
\<\rho((p(bb^*-1)+1)b^k)\xi\,|\,\rho(b^k)\xi\>
\nonumber\\ &=&
\<\rho(b^k)\xi\,|\,\rho(b^k)\xi\>-p^{k+1}\<\rho(b^k(1-bb^*))\xi\,|\,\rho(b^k)\xi\>
\nonumber\\
&=&
(1-p^{k+1})\<\rho(b^k)\xi\,|\,\rho(b^k)\xi\>\nonumber\\
\label{fakt} &=&(1-p^{k+1})\|\rho(b^k)\xi\|^2.
\eea
Thus $e_k:=\frac{\rho(b^k)\xi}{\|\rho(b^k)\xi\|}$ are normalised eigenvectors
of $\rho(1-bb^*)$ corresponding to the eigenvalue $p^k$. 
Note also that the vectors $e_k$ are orthogonal as they are eigenvectors
to different eigenvalues of the selfadjoint operator $1-\rho(bb^*)$.
Furthermore, it follows from the foregoing computation that
\beq\label{bek}
\rho(b)e_k=\frac{\rho(b^{k+1})\xi}{\|\rho(b^k)\xi\|}
=\frac{\rho(b^{k+1})\xi}{\|\rho(b^{k+1})\xi\|}
\frac{\|\rho(b^{k+1})\xi\|}{\|\rho(b^k)\xi\|}=\sqrt{1-p^{k+1}}\:e_{k+1}.
\eeq
On the other hand,
\[
\|\rho(b^*)\xi\|^2=\langle \rho(b^*)\xi\,|\,\rho(b^*)\xi\rangle
=\langle\rho(bb^*)\xi\,|\,\xi\rangle=0,
\]
so that $\rho(b^*)\xi=0$. Computing as in (\ref{fakt}),  we obtain
\beq\label{b*ek}
\rho(b^*)e_k=\frac{\rho(b^*b^k)\xi}{\|\rho(b^k)\xi\|}=\sqrt{1-p^k}\:e_{k-1},
\;\;\; k>0.
\eeq

Finally, it follows from $\cH'\inc\ker(\rho(1-aa^*))$ and (\ref{prel1})
that $\rho(a)$ is unitary on $\cH'$. Since it also commutes with $\rho(b)$
and $\rho(b^*)$,
 it belongs to the centre of the representation.
Hence, as the centre is always trivial in an irreducible representation, 
$\rho(a)$ has to be a multiple of the identity operator.
Therefore, the closed span of the orthonormal vectors $\{e_k\}_{k\geq0}$
is invariant under the whole algebra. Consequently, any irreducible 
$*$-representation on $\cH'$ is unitarily equivalent to one of the first
family of representations.
The second family is derived in the same way
exchanging the roles of $a$ and $b$, and $p$ and $q$.
\epf

Since, due to the disc-like relations (\ref{prel1})--(\ref{prel2}),
we have $\|\rho(a)\|\leq 1$ and $\|\rho(b)\|\leq 1$ for any bounded
$*$-representation of $\cO(S^3_{pq})$ \cite[Proposition IV.1(I)]{kl93}, 
  we can define
the \csa\ $C(S^3_{pq})$ of $\cO(S^3_{pq})$ as follows:
\bde\label{def}
$C(S^3_{pq})$ is the enveloping $C^*$-algebra of $\cO(S^3_{pq})$ for 
 the sup-norm
over all {\em bounded}
$*$-representations. 
\ede
Observe that
it follows from Theorem~\ref{irrep} that $\|\rho(a)\|=1=\|\rho(b)\|$ for any 
irreducible representation, so that  $\|a\|=1=\|b\|$. On the other hand,
using the unbounded representations of the quantum disc,
which can be obtained from the unbounded representations of the 
oscillator algebra given in \cite[Section 5.2.6]{ks97}, one obtains 
immediately two families of 
irreducible unbounded
$*$-representations $\rho_{1\theta\gamma}$  and $\rho_{2\theta\gamma}$ 
of $\cO(S^3_{pq})$.  
On a Hilbert space with an orthonormal basis $\{e_\mu\}_{\mu\in\Z}$,
they are given by the formulas
\beq
\rho_{1\theta\gamma}(a)e_\mu=\sqrt{1+q^{\mu+1}\gamma}\:e_{\mu+1},~~
\rho_{1\theta\gamma}(a^*)e_\mu=\sqrt{1+q^{\mu}\gamma}\:e_{\mu-1},~~
\gamma\in(q,1],
\eeq
\beq
\rho_{1\theta\gamma}(b)e_\mu=e^{i\theta}e_\mu,~~
\rho_{1\theta\gamma}(b^*)e_\mu=e^{-i\theta}e_\mu,~~\theta\in[0,2\pi),
\eeq
and analogous expressions for $\rho_{2\theta\gamma}$, with $a$ and $b$ exchanged
and $q$ replaced
by $p$.

To end this section, let us determine a vector space basis for ${\mathcal O}(S^3_{pq})$.

\begin{them}\label{basis} 
Let $\mu,\nu\in\Z$. Put $a_\mu=a^\mu$ if $\mu\geq 0$ and $a_\mu={a^*}^{|\mu|}$ if $\mu<0$.
Define $b_\nu$ in the same manner. Then
the elements 
\begin{equation}\label{S3-basis} 
\left\{ a_\mu(1-aa^*)^m(1-bb^*)^n b_\nu\;|\;\mu,\nu\in{\Bbb Z},\:m,n\in{\Bbb N},\:mn=0 \right\} 
\end{equation} 
form a vector space basis of ${\mathcal O}(S^3_{pq})$. Furthermore, 
${\mathcal O}(S^3_{pq})$ is faithfully embedded into $C(S^3_{pq})$. 
\end{them} 
\bpf 
At first we show that the elements (\ref{S3-basis}) linearly span ${\mathcal O}
(S^3_{pq})$. To this end, consider a word $W$ in $\{a,a^*,b,b^*\}$. Using relations 
(\ref{prel3}) we can rewrite $W$ as $W_aW_b$, where $W_a$ is a word in $\{a,a^*\}$ 
and $W_b$ is a word in $\{b,b^*\}$. Proceeding by induction on the length of 
$W_a$ 
and using the relation 
\begin{equation}
(1-aa^*)a=qa(1-aa^*),
\end{equation}
one can show that $W_a$ is a linear combination 
of polynomials of the form $a_\mu(1-aa^*)^m$. Likewise, using
\begin{equation}
(1-bb^*)b=pb(1-bb^*),
\end{equation}
$W_b$ is a linear combination 
of polynomials of the form $(1-bb^*)^n b_\nu$. 
Finally, since $a_\mu(1-aa^*)^m(1-bb^*)^n b_\nu=0$ unless 
$mn=0$ (see (\ref{prel4})), we can write any monomial $W$ as a linear 
combination of elements of (\ref{S3-basis}). Hence they span ${\mathcal O}(S^3_{pq})$.

Now, let us assume that $x$ is a finite linear combination of elements of (\ref{S3-basis}),
i.e.,
\beq
 x=\sum_{\mu,\nu\in\Z}\sum_{\scriptsize\begin{array}{l}m,n\in\N,\\ mn=0\end{array}}
 x_{\mu,m,n,\nu}a_\mu(1-aa^*)^m(1-bb^*)^n b_\nu. 
\eeq
In order to complete the proof of the theorem, we need to show the implication 
$$
\forall\;\mbox{bounded $*$-representations $\rho$}:\; \rho(x)=0
$$\[\label{imp}
\;\Downarrow\;
\]$$
\forall\;\mu,\nu\in\Z,\,m,n\in\N,\, mn=0:\;x_{\mu,m,n,\nu}=0.
$$
As only finitely of the coefficients $x_{\mu,m,n,\nu}$ are non-zero,
there exists $k\in\N$ such that $x_{\mu,m,n,\nu}=0$ for all $\mu<-k$.
Let us choose such a $k\in\N$ and assume
that $\rho_{1\theta}(x)=0$. Then it follows from (\ref{r1}) that
\beq
\rho_{1\theta}(x)e_k=\sum_{\scriptsize\begin{array}{l}\mu,\nu\in\Z,\\\nu\geq-k
\end{array}}\sum_{n\in\N}x_{\mu,0,n,\nu}e^{i\mu\theta}
(p^{k+\nu})^n\Lambda_{\nu k}e_{k+\nu}=0,
\eeq
where
\beq
\Lambda_{\nu k}=\left\{\begin{array}{cc}\sqrt{(1-p^k)\cdots(1-p^{k-\nu+1})}&-k\leq\nu<0\\[.3cm]
1&\nu=0\\[.3cm]
\sqrt{(1-p^{k+1})\cdots(1-p^{k+\nu})}&\nu>0
\end{array}\right..
\eeq
Since the vectors $e_i$ are linearly independent and 
$\Lambda_{\nu k}$'s are always
nonzero, for any $\nu\in\Z$ we have
\beq
\sum_{\mu\in\Z}\sum_{n\in\N}x_{\mu,0,n,\nu}e^{i\mu\theta}
(p^{k+\nu})^n=0.
\eeq
Next, as this equation is valid for all $\theta\in[0,2\pi)$, by the uniqueness of the 
Fourier coefficients, we can conclude that
\beq\label{pol}
\sum_{n\in\N}x_{\mu,0,n,\nu}(p^{k+\nu})^n=0
\eeq
for any $\mu,\nu\in\Z$. Since $p\not\in\{-1,0,1\}$ and there are infinitely many $k$'s
with the property $x_{\mu,m,n,\nu}=0$ for $\nu<-k$, the above polynomial in $p^{k+\nu}$
vanishes at infinitely many points, so that all its coefficients must be zero.
Thus we have shown that $x_{\mu,0,n,\nu}=0$ for all $\mu,\nu\in\Z,\,n\in\N$. 
The vanishing of the remaining coefficients $x_{\mu,m,0,\nu}$ can be proved using $\rho_{2\theta}$
instead of $\rho_{1\theta}$.

Finally, recall that if $x$ is annihilated by all representations, then its image
is zero in the enveloping $C^*$-algebra.
Therefore, the just proven implication (\ref{imp}) entails that
the elements $a^\mu(1-aa^*)^m(1-bb^*)^n b^\nu$ are linearly independent both as
elements of $C(S^3_{pq})$ and elements of $\cO(S^3_{pq})$. Hence (\ref{S3-basis})
is a basis of $\cO(S^3_{pq})$, and the canonical map
$\cO(S^3_{pq})\ra C(S^3_{pq})$ is injective.
\epf

\section{\boldmath$K$-theory of \boldmath$C(S^3_{pq})$}
\setcounter{equation}{0}

The purpose of this section is to show that the topological 
$K$-groups of $S^3_{pq}$ coincide
with the $K$-groups of the classical 3-sphere.
\bth\label{ks3}
The K-groups of the C*-algebra $C(S^3_{pq})$ are
$
K_0(C(S^3_{pq}))\cong\Z\cong K_1(C(S^3_{pq})).
$
\ethe
\bpf
Let $\cT$ denote the Toeplitz algebra, and $s$ its generating proper isometry
(the unilateral shift). 
The ideal of $\cT$ generated by $1-ss^*$ is isomorphic with the $C^*$-algebra 
$\cK$ of compact operators on a separable Hilbert space. We denote by $\pi$ 
the canonical surjection from $\cT\otimes\cT$ onto 
$(\cT\otimes\cT)/(\cK\otimes\cK)$. 
\ble\label{ttt} 
The $C^*$-algebras $C(S^3_{pq})$ and $(\cT\otimes\cT)/(\cK\otimes\cK)$ are 
isomorphic. 
\ele 
\bpf 
The universality of $C(S^3_{pq})$ for the relations 
(\ref{prel1})--(\ref{prel4}) and the representation formulas (\ref{r1})--(\ref{r2})
imply that there exists 
a $C^*$-algebra homomorphism $\alpha:C(S^3_{pq})\ra(\cT\otimes\cT)/
(\cK\otimes\cK)$ such that 
\bea 
&&\alpha(a) =  \pi(\rho_{2\theta}(a)\ot1)
=\pi\left(\sum_{n=0}^\infty\left(\sqrt{1-q^{n+1}}-
\sqrt{1-q^n}\right)s^{n+1}s^{*n}\otimes 1\right), \\ 
&&\alpha(b)  =  \pi(1\ot\rho_{1\theta}(b))
=\pi\left(\sum_{n=0}^\infty\left(\sqrt{1-p^{n+1}}-
\sqrt{1-p^n}\right)1\otimes s^{n+1}s^{*n}\right).  
\eea 
(We abuse notation by denoting with the same symbol representations of
$\cO(S^3_{pq})$ and $C(S^3_{pq})$.) Let us now construct the inverse of \ha. 
By (\ref{prel1}) we have $a^*a=1-q+qaa^*
\geq 1-q$, whence $a^*a$ is invertible. Therefore, so is $|a|=\sqrt{a^*a}$. 
Likewise, (\ref{prel2}) implies that $|b|$ is invertible. As both 
$a|a|^{-1}$ and $b|b|^{-1}$ are isometries, we can define a 
$C^*$-algebra homomorphism $\tilde{\beta}:\cT\otimes\cT\ra 
C(S^3_{pq})$ by 
\[ 
\tilde{\beta}(s\otimes 1)=  a|a|^{-1},\;\;\;
\tilde{\beta}(1\otimes s) = b|b|^{-1}. 
\] 
It follows from Theorem \ref{irrep} that $1$ is an isolated point of the 
spectrum of $1-aa^*$ with the corresponding spectral projection 
$1-a|a|^{-2}a^*$. Hence $1-a|a|^{-2}a^*$ belongs to the $C^*$-subalgebra 
of $C(S^3_{pq})$ generated by $1-aa^*$. Likewise, $1-b|b|^{-2}b^*$ belongs 
to the $C^*$-subalgebra of $C(S^3_{pq})$ generated by $1-bb^*$. 
Furthermore, the spectral projections are orthogonal  because the relation
(\ref{prel4}) implies that the eigenvectors $\psi_1$ and $\psi_2$, given by
$(1-aa^*)\psi_1=\psi_1$ and $(1-bb^*)\psi_2=\psi_2$, respectively, are
orthogonal:
\[
\<\psi_1\,|\,\psi_2\>
=
\<(1-aa^*)\psi_1\,|\,(1-bb^*)\psi_2\>
=
\<\psi_1\,|\,(1-aa^*)(1-bb^*)\psi_2\>
=0.
\]
 Thus  
$\tilde{\beta}((1-ss^*)\otimes(1-ss^*))=(1-a|a|^{-2}a^*)(1-b|b|^{-2}b^*)=0$. 
Since the smallest ideal of $\cT\otimes\cT$ containing $(1-ss^*)\otimes(1-ss^*)$ 
coincides with $\cK\otimes\cK$, we have $\tilde{\beta}(\cK\otimes\cK)=\{0\}$, 
and consequently $\tilde{\beta}$ induces a $C^*$-algebra homomorphism 
$\beta:(\cT\otimes\cT)/(\cK\otimes\cK)\ra C(S^3_{pq})$. It is straightforward
to verify on the generators that $\alpha\ci\beta=\id$:
\[
(\alpha\ci\beta)(\pi(s\ot1))
=
\ha(a|a|^{-1})
=
\pi(\rho_{2\theta}(a|a|^{-1})\ot1)
=
\pi(s\ot1).
\] 
(The case $1\ot s$ is analogous.) For the identity $\hb\ci\ha=\id$, note that
$\rho_{2\theta}$ and $\rho_{1\theta}$ are injective on the $C^*$-subalgebras
$C_a$ and $C_b$ generated by $a$ and $b$, respectively. Indeed, since $C(D_r)$
is the universal \csa\ for the relation $z^*z-rzz^*=1-r$, $r\in(0,1)$, we have
natural \csa\ epimorphisms $\pi_a:C(D_q)\ra C_a$ and $\pi_b:C(D_p)\ra C_b$.
On the other hand, $\rho_{2\theta}\ci\pi_a$ and $\rho_{1\theta}\ci\pi_b$
coincide with the faithful representation $\pi^I$ \cite[p.14]{kl93}, so that
$\rho_{2\theta}|_{C_a}$ and $\rho_{1\theta}|_{C_b}$ are injective. On the
other hand,
\bea
\rho_{2\theta}\llp(\hb\ci\ha)(a)-a\lrp
&=&
\rho_{2\theta}\llp \hb(\pi(\rho_{2\theta}(a)\ot1)-a\lrp
\nonumber\\ &=&
\sum_{n=0}^\infty(\sqrt{1-q^{n+1}}-
\sqrt{1-q^n})s^{n+1}s^{*n}\otimes 1
-\rho_{2\theta}(a)
\nonumber\\ &=&
0.
\eea
Similarly, $\rho_{1\theta}\llp(\hb\ci\ha)(b)-b\lrp=0$. Consequently, 
$\hb\ci\ha=\id$.
\epf\\
As shown in the foregoing lemma, there exists 
the following short exact sequence of $C^*$-algebras
\beq\label{ex}
0\lra
\cK\ot\cK\stackrel{j}{\lra}\cT\ot\cT\stackrel{\pi}{\lra}C(S^3_{pq})\lra 0,
\eeq
where $j$ is the inclusion map.
Applying the K\"unneth formula and remembering $K_0(\cK)\cong K_0(\cT)\cong \Z,
~K_1(\cK)\cong K_1(\cT)\cong 0$, reduces the six-term exact sequence
corresponding to (\ref{ex}) to
\beq\label{6ex}
0\lra K_1(C(S^3_{pq}))\stackrel{\partial}{\lra}\Z\stackrel{j_*}{\lra}
\Z\stackrel{\pi_*}{\lra}K_0(C(S^3_{pq}))\lra 0.
\eeq
Due to the exactness of this sequence, to end the proof it suffices to show that
$j_*=0$.
Note that $j_*$ goes from $K_0(\cK\ot\cK)$ to $K_0(\cT\ot\cT)$. Put 
$p=1-ss^*$. Then $p\ot p\in\cK\ot\cK$ is a minimal projection 
generating 
$K_0(\cK\ot\cK)$, and we have $j_*([p\ot p])=[p\ot p]\in K_0(\cT\ot\cT)$. 
In $K_0(\cT\ot\cT)$ we obviously have $[p\ot p]+[(1-p)\ot
p]=[1\ot p]$. On the other hand, 
$(s\ot p)(s\ot p)^*=(1-p)\ot p$ and $(s\ot p)^*(s\ot p)=1\ot p$, so that
$(1-p)\ot p$ and $1\ot p$ are Murray-von Neumann equivalent. Hence
$[1\ot p]=[(1-p)\ot p]$, which implies that $j_*=0$.
\epf

It was observed in \cite{hs1} that the $C^*$-algebra $C(SU_q(2))$ 
of the quantum 3-sphere of Woronowicz 
can be realized as the Cuntz-Krieger algebra corresponding to a 
certain finite directed graph. Similarly, $C(S^3_{pq})$  
can be identified with the algebra of a higher rank graph. We 
sketch this identification below, referring the interested reader 
to \cite{kp00} and \cite{rsy} for details on higher rank 
graphs. 
\bre\label{2g}\footnote{The following argument is due to Aidan Sims.}\em~
It is well-known that the Toeplitz algebra $\mathcal{T}$ is isomorphic 
to the $C^*$-algebra of the 1-graph $\Lambda$ with vertices $\Lambda^0 = 
\{v_1, v_2\}$ and oriented edges $\Lambda^1=\{e_1, e_2\}$ such that 
$r(e_1)=s(e_1)=r(e_2)=v_1$ and $s(e_2)=v_2$ (we use the edge-orientation 
convention of \cite{kp00}). The product $\Lambda\times\Lambda$ is a 
graph of rank 2 as in \cite[Definitions~2.1 and Definition~3.9]{rsy}, and  
the argument of \cite[Corollary~3.5 (iv)]{kp00} combined with 
\cite[Theorem~4.1]{rsy} shows that $C^*(\Lambda\times\Lambda) 
\cong C^*(\Lambda)\otimes C^*(\Lambda)\cong \mathcal{T} \otimes 
\mathcal{T}$. The Cartesian product 2-graph $\Lambda\times\Lambda$ has 
the 1-skeleton (see \cite[\S2]{rsy}) 

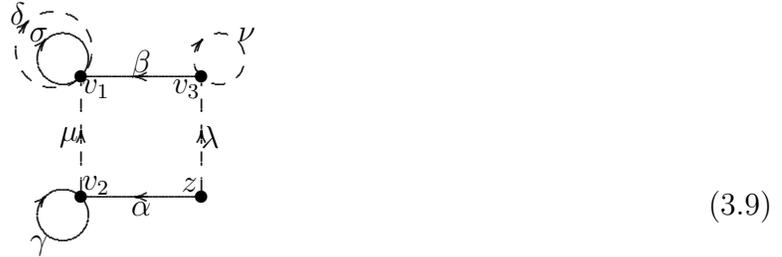
\vspace{3mm} 
\[
\beginpicture
\setcoordinatesystem units <.8cm,.8cm>
\put{$\bullet$} at 4 .3
\put{$\bullet$} at 6 .3
\put{$\bullet$} at 4 2.3
\put{$\bullet$} at 6 2.3
\circulararc 360 degrees from 3.4 -.3 center at 3.7 0
\arrow <0.15cm> [.25,.75] from 3.32  .21 to 3.41 .31
\arrow <0.15cm> [.25,.75] from 3.32  2.81 to 3.41 2.91
\circulararc 360 degrees from 3.4 2.9 center at 3.7 2.6
\plot 4.1 .3 5.9 .3 /
\plot 4.1 2.3 5.9 2.3 /
\arrow <0.15cm> [.25,.75] from 5.1  .3 to 4.9 .3
\arrow <0.15cm> [.25,.75] from 5.1  2.3 to 4.9 2.3
\arrow <0.15cm> [.25,.75] from 4  1.2 to 4 1.4
\arrow <0.15cm> [.25,.75] from 6  1.2 to 6 1.4
\setdashes
\plot 4 .4 4 2.2 /
\plot 6 .4 6 2.2 /
\circulararc 360 degrees from 3.1 3.2 center at 3.55 2.75
\circulararc 360 degrees from 6 2.9 center at 6.3 2.6
\arrow <0.15cm> [.25,.75] from 5.92  2.81 to 6.01 2.91
\arrow <0.15cm> [.25,.75] from 3.02  3.11 to 3.11 3.21
\setsolid
\put{$\delta$} at 2.95  3.35
\put{$\sigma$} at 3.3 3
\put{$\nu$} at 6.75 3
\put{$\gamma$} at 3.3 -.5
\put{$v_1$} at 4.25 2.1
\put{$v_2$} at 4.25 .5
\put{$v_3$} at 5.75 2.1
\put{$z$} at 5.8 .5
\put{$\lambda$} at 6.15 1.3
\put{$\mu$} at 3.8 1.3
\put{$\beta$} at 5 2.5
\put{$\alpha$} at 5 .1
\endpicture
\] 

\vspace{3mm}\noindent 
Here solid edges have degree (1,0) and dashed edges have degree (0,1). 
Let $\{s_\tau : \tau \in \Lambda\times\Lambda\}$ denote the universal 
Cuntz-Krieger $(\Lambda\times\Lambda)$-family. Now $\{z\}$ is a saturated 
hereditary subset of $(\Lambda \times \Lambda)^0$ in the sense of 
\cite[\S5]{rsy}. Writing $J$ for the closed ideal in $C^*(\Lambda\times\Lambda)$ 
generated by $S_z$, \cite[Theorem~5.2]{rsy} implies that $C^*(\Lambda\times
\Lambda)/J$ is isomorphic to the $C^*$-algebra of the quotient 2-graph, 
whose 1-skeleton is 

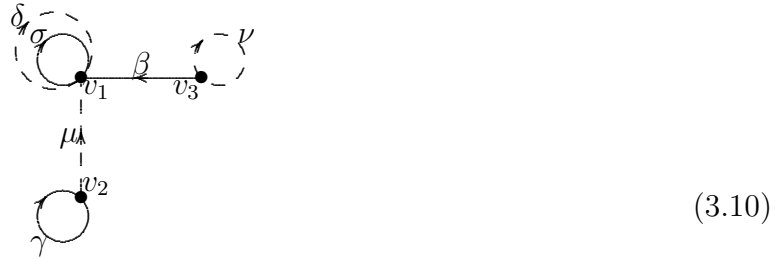
\vspace{3mm} 
\[
\beginpicture
\setcoordinatesystem units <.8cm,.8cm>
\put{$\bullet$} at 4 .3
\put{$\bullet$} at 4 2.3
\put{$\bullet$} at 6 2.3
\circulararc 360 degrees from 3.4 -.3 center at 3.7 0
\arrow <0.15cm> [.25,.75] from 3.32  .21 to 3.41 .31
\arrow <0.15cm> [.25,.75] from 3.32  2.81 to 3.41 2.91
\circulararc 360 degrees from 3.4 2.9 center at 3.7 2.6
\plot 4.1 2.3 5.9 2.3 /
\setdashes
\plot 4 .4 4 2.2 /
\circulararc 360 degrees from 3.1 3.2 center at 3.55 2.75
\circulararc 360 degrees from 6 2.9 center at 6.3 2.6
\arrow <0.15cm> [.25,.75] from 5.92  2.81 to 6.01 2.91
\arrow <0.15cm> [.25,.75] from 3.02  3.11 to 3.11 3.21
\setsolid
\arrow <0.15cm> [.25,.75] from 5.1  2.3 to 4.9 2.3
\arrow <0.15cm> [.25,.75] from 4  1.2 to 4 1.4
\put{$\delta$} at 2.95  3.35
\put{$\sigma$} at 3.3 3
\put{$\nu$} at 6.75 3
\put{$\gamma$} at 3.3 -.5
\put{$v_1$} at 4.25 2.1
\put{$v_2$} at 4.25 .5
\put{$v_3$} at 5.75 2.1
\put{$\mu$} at 3.8 1.3
\put{$\beta$} at 5 2.5
\endpicture
\] 

\vspace{3mm} \noindent
It is not difficult to verify that $I\cong\mathcal{K}\otimes\mathcal{K}$. 
Hence $(\mathcal{T}\otimes\mathcal{T})/(\mathcal{K}\otimes \mathcal{K})$
is the $C^*$-algebra of the above 2-graph. The former coincides with
$C(S^3_{pq})$ by Lemma~\ref{ttt}.
\ere

To end with, let us compare the geometry behind our computation and 
the corresponding calculations in \cite{mnw90} and \cite{m-k91a}. As can
be expected, in all three cases the $K$-groups are obtained from
the 6-term exact sequence of $K$-theory. Also in all three cases, they coincide
with their classical counterparts. However, the source of the 6-term exact sequence 
is each time different:
\cite[(0.2)]{mnw90}, \cite[p.355]{m-k91a}, (\ref{ex}). Considering the geometric
meaning of the employed \csa s (e.g., $\cK$ corresponds to $C_0(\R^2)$, $\cT$ to
$C(D)$, where $D$ is the unit disc in $\R^2$; see the Introduction in \cite{hms}),
we obtain the corresponding classical constructions:

\[
0\lra C_0(\R^2\times S^1)\lra C(S^3)\lra C(S^1)\lra 0\mbox{~~~(exact sequence),}
\]

\[
~~~~~~\begin{diagram}[height=10mm]
 C(S^3) &  \rTo           & C(D\times S^1) \\
 \dTo &  & \dTo \\
C(D\times S^1) & \rTo     & C(T^2) 
\end{diagram}
\mbox{~~~~~~(pullback diagram),}
\]

\[
0\lra C_0(\R^4)\lra C(D\times D)\lra C(S^3)\lra 0\mbox{~~~(exact sequence).}
\]

\noindent The first sequence means that removing $S^1$ from $S^3$ leaves a 
boundary-less solid torus. Think of $S^3=\{(c_1,c_2)\in D\times D\;|\;
|c_1|^2+|c_2|^2=1\}$ as a field of 2-tori over the internal points of
[0,1] bounded  by circles at the endpoints (e.g., take $|c_1|^2\in[0,1]$
as the interval parameter). Remove $S^1$ given by $c_1=0$ ($|c_2|=1$).
What remains is $S^1$ times a field of circles over the internal points of
(0,1] shrinking to a point at 1. The latter is an open disc (homeomorphic
with $\R^2$). The second sequence depicts gluing of two solid tori along their
boundaries, which is known as a Heegaard splitting of $S^3$. 
The corresponding six-term exact sequence is the Mayer-Vietoris sequence
of $K$-theory.
Finally, to
visualize the last sequence, recall that we can think of $S^3$ as the set
$\{ (z_1,z_2)\in D\times D\;|\; (1-|z_1|)(1-|z_2|)=0\}$ (see 
Proposition~\ref{classic} and divide by $(1+|z_1|)(1+|z_2|)$). 
Removing $S^3$ from $D\times D$ leaves all
points $(z_1,z_2)\in D\times D$ such that $(1-|z_1|)(1-|z_2|)\neq 0$,
which is precisely the Cartesian product of two open discs (homeomorphic with
$\R^4$).

\section{Hopf-Galois aspects of \boldmath$\cO(S^2_{pq})\subset \cO(S^3_{pq})$}
\setcounter{equation}{0}

Our goal now is to prove that the extension $\cO(S^2_{pq})\subset \cO(S^3_{pq})$
is Hopf-Galois, relatively projective, and non-cleft. We begin with the following:
\ble\label{galois}
The locally trivial $\cO(U(1))$-extension $\cO(S^2_{pq})\subset \cO(S^3_{pq})$
is Hopf-Galois.
\ele
\bpf
Note first that, since $\cO(U(1))$ is cosemisimple, the injectivity of
the canonical map (see Subsection~1.3)
\[
can:{\cal O}(S^3_{pq})\ot_B{\cal O}(S^3_{pq}) 
\lra {\cal O}(S^3_{pq})\ot {\cal O}(U(1))
\]
follows from its surjectivity. Indeed, since there is a Haar functional $f_H$ on $\cO(U(1))$, 
we get the total integral
of Doi by composing it with the unit map: $j:=\eta\ci f_H:\cO(U(1))\ra\cO(S^3_{pq})$.
Therefore, we can apply Remark~3.3 and
Theorem~I of \cite{s-hj90}. On the other hand, as explained in Subsection~1.3, 
to prove the surjectivity of $can$, it suffices to show that both $1\ot u$ and $1\ot u^*$
are in its image. (Here $u$ and $u^*$ are the generators of $\cO(U(1))$.) Taking advantage of
(\ref{dpa}), the relations (\ref{prel1}), (\ref{prel3}), and the key relation (\ref{prel4}), 
we obtain
\bea\label{can1}
can(a^*a\ot_B u+qbb^*(1-aa^*)\ot_B u)
&=&(a^*a+qbb^*-qaa^*bb^*)\ot u\nonumber\\
&=&(qaa^*+1-q+qbb^*-qaa^*bb^*)\ot u\nonumber\\
&=&1\ot u.
\eea
 Analogously, using (\ref{prel2})
instead of (\ref{prel1}), we obtain
\bea\label{can2}
can(b^*b\ot_B u^*+paa^*(1-bb^*)\ot_B u^*)
&=&(b^*b+paa^*-pbb^*aa^*)\ot u^*\nonumber\\
&=&(pbb^*+1-p+paa^*-pbb^*aa^*)\ot u^*\nonumber\\
&=& 1\ot u^*.
\eea
Thus we have shown that
$can$ is bijective, as needed.
\epf

Our next step is to construct a strong connection $\ell$ (see Subsection~1.3).
\ble\label{rp}
The Hopf-Galois $\cO(U(1))$-extension $\cO(S^2_{pq})\inc\cO(S^3_{pq})$ is relatively
projective.
\ele
\bpf
Since, due to
\cite{bh}, the \hge\ is relatively projective if and only if there exists a strong 
connection, we prove the assertion of the lemma by constructing a strong connection $\ell$.
We define a linear map $\ell:\cO(U(1))\ra 
\cO(S^3_{pq})\ot\cO(S^3_{pq})$, $\ell(h)=h^{[1]}\ot h^{[2]}$ (summation understood),
by giving its values on the basis elements
$u^\mu$, $\mu\in \Z$:
\[
\ell(1)=1\ot 1,
\]
\[
\ell(u)=a^*\ot a + qb(1-aa^*)\ot b^*,\;\;\;
\ell(u^*)= b^*\ot b +pa(1-bb^*)\ot a^*,
\]
\[
\ell(u^\mu)=u^{[1]}\ell(u^{\mu-1})u^{[2]},\;\;\;
\ell(u^{*\mu})=u^{*[1]}\ell(u^{*(\mu-1)})u^{*[2]},\;\;\;\mu>0.
\]
Denote by $\widetilde{can}$ the lifting $(m\ot\id)\ci(\id\ot\dr)$ of the
canonical map $can$. (Here $m$ stands for the multiplication map.)
To show that $\tc\ci\ell=1\ot \id$, we first note that (\ref{can1})  and (\ref{can2})
entail that $(\tc\ci \ell)(u)=1\ot u$ and $(\tc\ci \ell)(u^*)=1\ot u^*$.
Assume now that $(\tc\ci \ell)(u^k)=1\ot u^k$. Then
\bea
(\tc\ci \ell)(u^{k+1})&=& \tc(u^{[1]}\ell(u^k)u^{[2]})\nonumber\\
&=&u^{[1]}((\tc\ci \ell)(u^k))\Delta_R(u^{[2]})\nonumber\\
&=&u^{[1]}(1\ot u^k)\Delta_R(u^{[2]})\nonumber\\
&=&1\ot u^{k+1}.
\eea
The case $(\tc\ci \ell)({u^*}^k)=1\ot{u^*}^k$ can be handled in the same way.
Therefore, it follows by induction that $\tc\ci\ell=1\ot \id$.
Now, let us consider the following diagram (cf.\ \cite[(1.25)]{dgh01}:
\begin{diagram}[height=10mm]
         &           &\cO(S^3_{pq})\ot\cO(S^3_{pq})                &          & \\
         &\ruTo{\ell}&\dTo_{\pi_B}                                 &\rdTo{\tc}& \\
\cO(U(1))&\rTo{T}    &\cO(S^3_{pq})\ot_{\cO(S^2_{pq})}\cO(S^3_{pq})&\rTo{can} &
\cO(S^3_{pq})\ot\cO(U(1)) . 
\end{diagram}
The right triangle part of the diagram commutes by construction, and we have already
shown that the big triangle commutes. Thus the commutativity of the left triangle
follows from the injectivity of $can$. This means that
$\ell$ is  a lifting of the translation map. It is by construction unital, so that
it remains to show its bicolinearity, which is again done inductively.
First, it is immediate to see the equality 
$(\ell\ot\id)\ci\Delta=(\id\ot\Delta_R)\ci \ell$ on the generator $u$.
 We can write this as
\beq\label{u12}
u^{[1]}\ot u^{[2]}\ot u=u^{[1]}\ot\Delta_R(u^{[2]}).
\eeq
Next, assume that $((\ell\ot\id)\ci\Delta)(u^k)=((\id\ot\Delta_R)\ci\ell)(u^k)$
for some $k>0$. Then 
\[
((\ell\ot\id)\ci\Delta)(u^{k+1})
=
\ell(u^{k+1})\ot u^{k+1}
=
u^{[1]}\ell(u^k)u^{[2]}\ot u^{k+1}.
\]
On the other hand, taking advantage of the inductive assumption, we obtain
\bea
((\id\ot\dr)\ci\ell)(u^{k+1})
&=&(\id\ot\Delta_R)(u^{[1]}\ell(u^k)u^{[2]})\nonumber\\
&=&u^{[1]}\llp(\id\ot\Delta_R)(\ell(u^k))\lrp\Delta_R(u^{[2]})\nonumber\\
&=&u^{[1]}(\ell\ot\id)(\Delta(u^k))\Delta_R(u^{[2]})\nonumber\\
&=&u^{[1]}(\ell(u^k)\ot u^k)\Delta_R(u^{[2]})\nonumber\\
&=&u^{[1]}\ell(u^k)u^{[2]}\ot u^{k+1}.
\eea
A similar argument can be made with $u^*$ in place of $u$.
This proves the right colinearity of $\ell$ due to the fact that
$\{u^\mu\}_{\mu\in \Z}$ is a basis of $\cO(U(1))$.
The proof of the left colinearity is fully analogous, now using
\beq
u\ot u^{[1]}\ot u^{[2]}=\Delta_L(u^{[1]})\ot u^{[2]}
\eeq
and its $*$-version. Thus we have shown  that $\ell$ is a strong connection. 
\epf

\bre\em
Galois coactions are algebraic incarnations of principal actions in 
classical geometry. There is an attempt to find the corresponding principality
or Galois condition for \csa s \cite{e-da00}. In our situation, Definition~2.4
of \cite{e-da00} defining principal coactions on \csa s reduces to requiring that
$C(S^3_{pq})\dr(C(S^3_{pq}))$ be norm dense in $C(S^3_{pq})\ot C(U(1))$. Since
the powers of $u$ span a dense subspace of $C(U(1))$, it suffices to note
that $1\ot u^\mu\in C(S^3_{pq})\dr(C(S^3_{pq}))$. This follows from
the proof of the foregoing theorem, so that 
the coaction $\Delta_R:C(S^3_{pq})\ra C(S^3_{pq})\ot C(U(1))$ 
is principal in the sense of \cite{e-da00}. 
\ere

We now provide an explicit formula for $\ell$.
Put $x=a^*\ot a$ and $y=qb(1-aa^*)\ot b^*$, and consider them 
as elements of $\cO(S^3_{pq})^{op}\ot\cO(S^3_{pq})$, where 
the superscript $op$ indicates the opposite algebra. Then the formula 
for $\ell(u^k)$ reads  
$
l(u^k)=(x+y)^k. 
$ 
Due to (\ref{prel1}),
$
yx=qxy, 
$ so that
we can use the formula 
\beq\label{lxy} 
(x+y)^n=\sum_{k=0}^n\mbox{\scriptsize$\left(\!\ba{cc}n\\k\ea\!\right)_q$}x^ky^{n-k},
\;\;\; 
\left(\!\ba{cc}n\\k\ea\!\right)_q
=\frac{(q-1)\cdots (q^n-1)}{(q-1)\cdots (q^k-1) (q-1)\cdots (q^{n-k}-1)}. 
\eeq 
Consequently, we obtain
\beq\label{lex} 
\ell(u^n)=\sum_{k=0}^n\mbox{\scriptsize$\left(\!\ba{cc}n\\k\ea\!\right)_q$}
q^{n-k}(1-aa^*)^{n-k}a^{*^k} 
b^{n-k}\ot a^kb^{*^{n-k}}. 
\eeq 
The formula for $\ell(u^{*^n})$ can be derived exchanging the roles of
 $a$ and $b$,  and $q$ and $p$.
Notice also that, since $\tc\ci\ell=\id$, we have
\[
m\ci\ell=(\id\ot\he)\ci\tc\ci\ell=\he.
\] 
Thus (\ref{lex}) entails the following identity in $\cO(S^3_{pq})$:
\[
\sum_{k=0}^n\mbox{\scriptsize$\left(\!\ba{cc}n\\k\ea\!\right)_q$}
q^{n-k}(1-aa^*)^{n-k}a^{*^k} 
b^{n-k}a^kb^{*^{n-k}}=1.
\]
A similar identity follows from an explicit formula for $\ell(u^{*k})$.


Next, consider the 1-dimensional corepresentations of $\cO((U(1)))$,
$\rho_\mu(1)=1\ot u^{-\mu}$, $\mu\in\Z$. We can identify 
$\hom_{\rho_\mu}(\C,\cO(S^3_{pq}))$ with
$\cO(S^3_{pq})_\mu:=\{p\in\cO(S^3_{pq})\;|\;\dr(p)=p\ot u^{-\mu}\}$.
Since the powers of $u$ form a basis of $\cO((U(1)))$, we have the direct
sum decomposition
$\cO(S^3_{pq})=\bigoplus_{\mu\in\Z} \cO(S^3_{pq})_{\mu}$
as $\cO(S^2_{pq})$-bimodules.
According to the general result referred to in Subsection~1.3,
the strong connection determines idempotent matrices $E_\mu$ of the associated 
modules ($\cO(S^3_{pq})_\mu\cong \cO(S^2_{pq})^{{\rm size}E_\mu}E_\mu$
as $\cO(S^2_{pq})$-modules). 
Using   formula (\ref{lex}), one finds explicitly  the idempotents 
$E_{-n}=R_{-n}^TL_{-n}$, $n\in\N$, where
\beq\label{e-} 
R_{-n}= 
({b^*}^n, a{b^*}^{n-1},\ldots,a^n),
\eeq
\beq
L_{-n}=\left(
\mbox{\scriptsize$\left(\!\ba{cc}n\\0\ea\!\right)_q$}q^n(1-aa^*)^nb^n\,\bc\,
\mbox{\scriptsize$\left(\!\ba{cc}n\\1\ea\!\right)_q$}q^{n-1}(1-aa^*)^{n-1}a^*b^{n-1}
\,\bc\,\cdots\,\bc\,
\mbox{\scriptsize$\left(\!\ba{cc}n\\n\ea\!\right)_q$}{a^*}^n\right).
\eeq
The idempotents $E_n$ are given by a similar formula, with $a$ and $b$ exchanged and  
$q$ replaced by $p$. 
For $\mu=-1$, we have
\[\label{e1}
E_{-1}:=\left(\ba{c} a\\b^*\ea\right)
\left(\ba{cc} a^*&qb(1-aa^*)\ea\right)
=\left(\ba{cc}
aa^*&qa(1-aa^*)b\\
a^*b^*&q(1-aa^*)b^*b\\
\ea\right).
\] 

Now we want to prove
 that the \hge\ $\cO(S^2_{pq})\inc\cO(S^3_{pq})$ is not cleft. We do it by showing
that the $K_0$-class of the idempotent $E_{-1}$ is not trivial. This, in turn,
we prove by computing an appropriate invariant of $K$-theory. The invariant has
a very simple form, namely it is the Chern-Connes pairing of a trace 
(0-cyclic cocycle) with $E_{-1}$. It is known
\cite{hmsf} that
there exists a trace on $\cO(S^2_{pq})$ given by
\begin{equation} \label{tr}
{\rm tr}(f):={\rm Tr}(\rho_2(f)-\rho_1(f)).
\end{equation} 
Here $\rho_1,\rho_2:\O(S^2_{pq})\rightarrow\B(\H)$ are the two 
infinite-dimensional 
representations given by (1.21)--(1.22), and {\rm Tr} is the operator 
trace.
The pairing of tr and $E_\mu$ has been computed in
\cite{hmsf} for any $\mu\in\Z$. Since the special case $\mu=-1$ computation
is straightforward, we enclose it here for the convenience of the reader.
\ble\label{CC}
Let $\langle \;,\;\rangle$ denote the pairing between the cyclic cohomology
$HC^{even}(\cO(S^2_{pq}))$ and $K_0(\cO(S^2_{pq}))$. Then
$
\< \tr, [E_{-1}] \>  = -1.
$
\ele
\bpf
Using (\ref{e1}) and (\ref{prel1})--(\ref{prel4}),
we have 
\bea
\< \tr, [E_{-1}] \>  
&=&
\tr(\Tr_{M_2}(E_{-1}))
\nonumber\\ &=&
\tr(aa^*+q(1-aa^*)b^*b)
\nonumber\\ &=&
\tr(aa^*+q(1-aa^*)(p(bb^*-1)+1))
\nonumber\\ &=&
\tr(q+(1-q)aa^*).
\eea
 Taking advantage of (\ref{coin}), this can be expressed in terms of $f_0$ and
$f_1$. Using again the commutation relations, we get
\[
aa^*=1-bb^*+ab(ab)^*=1-f_0+f_1f_{-1}\;\;\;\mbox{(injection $\iota$ suppressed).}
\]
It is immediate from (\ref{tr}) that $\tr(1)=0$,  and it follows from 
(\ref{ro1})--(\ref{ro2}) 
that $\tr(f_0-f_1f_{-1})=\frac{1}{1-q}$.
This yields
\[
\tr(q+(1-q)aa^*)=(1-q)\tr(aa^*)=(q-1)\tr(f_0-f_1f_{-1})=-1,
\]
as needed.
\epf

Since every free module can be represented in $K_0$ by the identity
matrix, the pairing between $\tr$ and the $K_0$-class of any free 
module always yields zero.
Thus the left module $\cO(S^3_{pq})_{-1}\cong \cO(S^2_{pq})^2E_{-1}$ is not (stably) free.
Now, reasoning as in \cite[Section~4]{hm99}, we can conclude that the $\cO(U(1))$-extension
$\cO(S^2_{pq})_{-1}\inc \cO(S^3_{pq})$ is non-cleft.
Combining this with Lemma~\ref{galois} and Lemma~\ref{rp}, we obtain:
\bth\label{noncl}
The locally trivial  $\cO(U(1))$-extension $\cO(S^2_{pq})\inc\cO(S^3_{pq})$ is 
a relatively projective non-cleft Hopf-Galois extension.
\ethe

\bre\em
Note that in \cite{cm} it was only shown that $\cO(S^3_{pq})$ is 
not isomorphic to the tensor product $\cO(S^2_{pq})\otimes \cO(U(1))$
(non-triviality). Here we prove that $\cO(S^3_{pq})$ is not a crossed product
$\cO(S^2_{pq})\rtimes \cO(U(1))$. (For a discussion concerning non-trivial 
versus
non-cleft, see the end of Section~4 in \cite{dhs99}.)
\ere

{\bf Acknowledgments.}
The authors are indebted to A.\ Sims and  A.\ Sitarz for their help with Remark~\ref{2g}
and Proposition~\ref{classic}, respectively, and
to T.\ Brzezi\'nski and P.\ Schauenburg for their help with the proofreading.
P.M.H.\ has been supported by a Marie Curie Fellowship of the European
Community under the contract number HPMF-CT-2000-00523. P.M.H.\ is also
grateful to the Naturwissenschaftlich-Theoretisches Zentrum der 
Universit\"at Leipzig and the University of Newcastle, Australia,
for their hospitality and financing his visit to Leipzig and Newcastle, respectively.
R.M.\ has been supported by the Deutsche Forschungsgemeinschaft.
R.M.\ is also grateful for hospitality and financial support to the 
Universit\"at M\"unchen.  
R.M.\ and W.S.\ thank 
the European Commission for covering the expenses of their stay in Warsaw
during the Banach Centre school/conference ``Noncommutative Geometry and Quantum
Groups.''
W.S.\ also would like to thank the Research Management Committee of the University
of Newcastle and the Max-Planck-Institut f\"ur Mathematik Bonn for their support.
All three authors are grateful to the Mathematisches Forschungsinstitut
Oberwolfach for support via its Research in Pairs programme.

\end{document}